\documentclass[a4paper,reqno]{amsart}
\usepackage[latin1]{inputenc}
\usepackage{upref}
\pdfoutput=1
\usepackage{graphicx}
\usepackage{hyperref}

\newtheorem{theorem}{Theorem}[section]
\newtheorem{lemma}[theorem]{Lemma}
\newtheorem{corollary}[theorem]{Corollary}
\newtheorem{hypothesis}[theorem]{Hypothesis}
\theoremstyle{definition}
\newtheorem{definition}[theorem]{Definition}

\theoremstyle{remark}
\newtheorem{remark}[theorem]{Remark}
\allowdisplaybreaks
\numberwithin{equation}{section}

\newcommand{\dott}{\,\cdot\,}
\newcommand{\seq}[1]{\left\{#1\right\}}
\newcommand{\norm}[1]{\left\|#1\right\|}
\newcommand{\dx}{\,dx}
\newcommand{\dy}{\,dy}
\newcommand{\abs}[1]{\left|#1\right|}
\newcommand{\Dx}{\Delta x}
\newcommand{\Dt}{\Delta t}
\newcommand{\Dp}{D_+}
\newcommand{\Dm}{D_-}
\newcommand{\Dpm}{D_\pm}

\newcommand{\eps}{\varepsilon}
\newcommand{\R}{\mathbb{R}}
\newcommand{\Z}{\mathbb{Z}}

\newcommand{\Pt}{\Pi_T}

\newcommand{\Dtp}{D^t_+}

\newcommand{\weakto}{\rightharpoonup}

\newcommand{\test}{\varphi}

\newcommand{\sgn}[1]{\mathrm{sign}\left(#1\right)}

\newcommand{\cha}[1]{\mathbf{1}_{#1}}
\newcommand{\tc}{\tilde{c}}

\newcommand{\wlim}[1]{\overline{#1}}
\newcommand{\cj}[1]{c_{j#1 1/2}}
\newcommand{\uj}[1]{u_{j#1 1/2}}
\newcommand{\udx}{u_{\Dx}}
\newcommand{\Rdx}{R_{\Dx}}
\newcommand{\Sdx}{S_{\Dx}}
\newcommand{\Le}{\mathcal{L}_{\eps}}
\newcommand{\Sobolev}[2]{W_{\mathrm{loc}}^{#1,#2}}
\newcommand{\Hneg}{H_{\mathrm{loc}}^{-1}}

\newcommand{\Lenloc}{L_{\mathrm{loc}}^{1}}
\newcommand{\Mloc}{\mathcal{M}_{\mathrm{loc}}}
\newcommand{\Div}{\mathrm{div}}
\newcommand{\Curl}{\mathrm{curl}}
\newcommand{\weakstar}{\overset{\star}\rightharpoonup}
\newcommand{\Fjmh}{F_{j-1/2}}
\newcommand{\Fjph}{F_{j+1/2}}
\newcommand{\ujmh}{u_{j-1/2}}
\newcommand{\ujph}{u_{j+1/2}}
\newcommand{\el}[1]{L^{#1}(\R)}
\newcommand{\cjph}{c_{j+1/2}}
\newcommand{\cjmh}{c_{j-1/2}}

\newcommand{\loc}{\mathrm{loc}}
\newcommand{\weak}{\rightharpoonup}

\newcommand{\Set}[1]{\left\{#1\right\}}
\newcommand{\supp}{\mathrm{supp}}

\title[Finite differences for  a nonlinear  wave equation]{A convergent finite difference method 
for a nonlinear variational wave equation}
\author[Holden]{H.~Holden}
\address[Holden]{\newline
    Department of Mathematical Sciences,
    Norwegian University of Science and Technology,
    NO--7491 Trondheim, Norway,
and \newline
  Centre of Mathematics for Applications, 
University of Oslo,
  P.O.\ Box 1053, Blindern,
  NO--0316 Oslo, Norway  }
\email[]{holden@math.ntnu.no}
\urladdr{www.math.ntnu.no/\~{}holden/}
\author[Karlsen]{K.~H.~Karlsen}
\address[Karlsen]{\newline
   Centre of Mathematics for Applications, 
 University of Oslo,
  P.O.\ Box 1053, Blindern,
  NO--0316 Oslo, Norway}
\email[]{kennethk@math.uio.no}
\urladdr{www.math.uio.no/\~{}kennethk/}

\author[Risebro]{N.~H.~Risebro}
\address[Risebro]{\newline
  Centre of Mathematics for Applications, 
University of Oslo,
  P.O.\ Box 1053, Blindern,
  NO--0316 Oslo, Norway}
\email[]{nilshr@math.uio.no}
\urladdr{www.math.uio.no/\~{}nilshr/}

\date{\today}
\thanks{This research was supported in part by the Research Council of Norway. 
KHK has been supported in part by an Outstanding Young Investigators 
Award from the Research Council of Norway.}

\begin{document}

\begin{abstract}  
We establish rigorously convergence of a semi-discrete upwind scheme for the 
nonlinear variational wave equation $u_{tt} - c(u)(c(u) u_{x})_x = 0$ with 
$u|_{t=0}=u_0$ and $u_t|_{t=0}=v_0$.  Introducing  Riemann invariants 
$R=u_t+c u_x$ and $S=u_t-c u_x$, the variational wave equation is equivalent to 
$R_t-c R_x=\tilde c (R^2-S^2)$ and $S_t+c S_x=-\tilde c (R^2-S^2)$ with  $\tilde c=c'/(4c)$. 
An upwind scheme is defined for this system. We assume that the the speed $c$ is 
positive, increasing and both $c$ and its derivative are 
bounded away from zero and that $R|_{t=0}, S|_{t=0}\in L^1(\R)\cap L^3(\R)$ 
are nonpositive. The numerical scheme is illustrated on several examples. 
\end{abstract}

\maketitle
\section{Introduction}\label{sec:intro}
In this paper we consider the nonlinear variational wave equation
\begin{equation}
  \label{eq:waveeq}
  \begin{gathered}
    \frac{\partial^2 u}{\partial t^2} - c(u)\frac{\partial}{\partial
      x}\left(c(u)\frac{\partial u}{\partial x}\right) = 0,\\
    u(x,0)=u_0(x),\quad \frac{\partial u}{\partial t}(x,0)=v_0(x),
  \end{gathered}
\end{equation}
in the strip $(x,t)\in \Pt=\R\times [0,T]$. 

The equation, which can be derived as the Euler--Lagrange equation for
the variational principle $\delta
\iint(\psi_{t}^2-c^2(\psi)\psi_{x}^2)\, dxdt=0$, can be used to model
liquid crystals, see \cite{Saxton1989,HunterSaxton1991,
GlasseyHunterZheng:sing}.  Consider namely a nematic liquid
crystal in the regime where inertial effects dominate. In that case
the liquid crystal can be described by the director field
$n=n(x,y,z,t)\in\R^3$ with $\norm{n}=1$ that describes the direction
of the elongated molecules that constitute the liquid crystal. Its
potential energy density is described by the Oseen--Franck functional
\begin{equation*}
W(n,\nabla n)= \alpha \abs{n\times 
(\nabla\times n)}^2+\beta (\nabla\cdot n)^2
+\gamma (n\cdot\nabla n)^2,
\end{equation*}
where the constants $\alpha$, $\beta$, and $\gamma$ describe the
liquid crystal. The dynamics of the director field is given by a least
action principle
\begin{equation}\label{eq:leastaction}
\frac{\delta}{\delta u}\int (n_t\cdot n_t- W(n,\nabla n))\, dxdydz dt =0.
\end{equation}
Consider the class of planar deformations given by   
\begin{equation}
n=\cos(u(x,t)) \mathbf{i}+ \sin(u(x,t)) \mathbf{j}
\end{equation}
where $ \mathbf{i}$ and $\mathbf{j}$ are unit vectors in the $x$ and
$y$ direction, respectively.  In that case the least action principle
\eqref{eq:leastaction} reduces to \eqref{eq:waveeq} with
\begin{equation*}
c^2(u)=\alpha \cos^2 u+\beta \sin^2 u.
\end{equation*}
We here analyze \eqref{eq:waveeq} with more restrictive assumptions on
$c$, as is done in the mathematical literature, namely that $c$ is
positive, strictly increasing and bounded away from zero.  We note that 
\eqref{eq:waveeq} is closely related to the Hunter--Saxton
equation, which is obtained by a further asymptotic expansion of
\eqref{eq:waveeq}, see \cite{HunterSaxton1991}.  


While short-term existence of regular solutions follows by the Kato
method, it is clear that the solution in general develops
singularities in finite time, even from smooth initial data, see
\cite{GlasseyHunterZheng,GlasseyHunterZheng:sing}.

In a series of papers, Zhang and Zheng \cite{ZZ:2001a,
  ZZ:rarefactive,ZZ:2003:weakvariational,ZZ:2004,ZZ:2005} have
analyzed \eqref{eq:waveeq} carefully, using the method of Young
measures.  From their many results we quote the recent one
\cite[Thm.~1.1]{ZZ:2004} where they show existence of a global weak
solution for initial data $u_0\in H^1(\R)$ and $v_0\in L^2(\R)$.  The
function $c$ is assumed to be smooth, bounded, positive with
derivative that is non-negative and strictly positive on the initial
data $u_0$. Their results, and also the relationship to the
Hunter--Saxton equation are surveyed in \cite{ZZ:2005}. The uniqueness
question is open.

Another approach to the study of \eqref{eq:waveeq} was recently taken
by Bressan and Zheng \cite{BressanZheng2005}.  Instead of following
the approach based on Young measures, they rewrite the equation in new
variables where singularities disappear.  They show that for $u_0$
absolutely continuous with $u_{0,x},v_0\in L^2(\R)$ the Cauchy-problem
\eqref{eq:waveeq} allows a global weak solution with the following
properties: The solution $u$ is locally H\"older continuous with
exponent $\frac12$, and the map $t\mapsto u(t,\dott)$ continuously
differentiable with values in $L^p_{\mathrm{loc}}(\R)$ for $1\le p<2$.
Further properties are obtained, in particular, it is shown that the
associated energy, treated as a measure, is conserved in time.

Up to now, little has been known about the behavior of numerical
schemes for the equation \eqref{eq:waveeq}.  Except for some numerical
computations in \cite{GlasseyHunterZheng:sing}, there are, to the best
of our knowledge, no rigorous results regarding any numerical method
for \eqref{eq:waveeq}, and the main purpose of this paper is to remedy
this situation. Here we introduce a semi-discrete upwind scheme for
the initial-value problem \eqref{eq:waveeq}, i.e., a finite difference
approximation of the spatial variation, keeping the continuous
temporal variation.  For this scheme we show convergence to a weak
solution of \eqref{eq:waveeq}, and thus this proof offers a different
existence proof compared with the others. In addition it provides a
constructive approach to the initial-value problem in the sense that
the difference scheme supplies a numerical tool to compute the
solution, see Section~\ref{sec:numer}.  Indeed, we show how the
difference scheme performs, both on examples where the scheme is
proved to converge and otherwise.  A similar analysis has been applied
to the Hunter--Saxton equation, see~\cite{HoldenKarlsenRisebro}.

\medskip We now turn to a more detailed and technical discussion. Weak
solutions are defined as follows. 
\begin{definition} \label{def:weaksol} 
  Let $\Pt$ denote the set $\R\times [0,T)$, $T>0$.
  By a weak solution $u$ of \eqref{eq:waveeq} we mean a function
  $u\in L^\infty([0,T]; W^{1,p}(\R))\cap C(\Pt)$, $u_t \in
  L^\infty([0,T]; L^{p}(\R))$, for all $p\in[1,3+q]$,
  where $q$ is some fixed positive constant $q>0$, such that 
  \begin{equation}
    \label{eq:weaksol}
    \iint\limits_{\Pt}\big( \partial_t \test\, \partial_t u - c^2(u)\partial_x
    \test \,\partial_x u - c'(u)c(u)\test \left(\partial_x u\right)^2\big) 
    \,dxdt = 0
  \end{equation}
  for all test functions $\test\in C^\infty_0(\Pt)$. The initial values
  are taken in the sense that $u(\dott,t)\to u_0$ in $C([0,T];L^2(\R))$
  as $t\to 0+$, and $\partial_t u(\dott,t)\to v_0$ as a distribution in
  $\Pt$ when $t\to 0+$.
\end{definition} 
A common approach to \eqref{eq:waveeq} is first to re-write the
equation in terms of Riemann invariants. To that end we define
\begin{equation*}
R=\frac{\partial u}{\partial t} + c(u)\frac{\partial u}{\partial
  x},\quad S=\frac{\partial u}{\partial t} - c(u)\frac{\partial u}{\partial x},
\end{equation*}
and the auxiliary function
\begin{equation*}
\tc(u)=\frac{c'(u)}{4c(u)}.
\end{equation*}
Then the wave equation \eqref{eq:waveeq} is formally equivalent to the
$3\times 3$ system 
\begin{align}
  \label{eq:wavesys}
  &\begin{aligned}
    R_t - c(u)R_x &= \tc(u)\left(R^2-S^2\right),\\
    S_t + c(u)S_x &= -\tc(u)\left(R^2-S^2\right),\\
    u_x &= \frac{1}{2c(u)}(R-S), 
  \end{aligned}
  \\
 & R_0=R|_{t=0}= v_0+c(u_0)u_0', \qquad S_0=S|_{t=0}= v_0-c(u_0)u_0'.
  \label{eq:R0S0def}
\end{align}
Clearly, we also have
\begin{equation}\label{eq:ux}
    u_t=\frac12(R+S).
\end{equation}
In order to make this well defined, we use the boundary condition
\begin{equation*}
\lim_{x\to-\infty} u(x,t)=0.
\end{equation*}
The equations for $R, S$  can also be written on conservative form, viz.
\begin{equation}
  \label{eq:wavesysCONS}
  \begin{aligned}
    R_t - \big(c(u)R\big)_x &= -\tc(u)\left(R-S\right)^2,\\
    S_t + \big(c(u)S\big)_x &= -\tc(u)\left(R-S\right)^2.
  \end{aligned}
\end{equation}
Throughout the paper we will assume that  $c$ is a Lipschitz continuous function such that
\begin{equation}
0<C_1 \le c(u) \le C_2,\quad\text{and}\quad 0\le c'(u) \le M_1.\label{c_est}
\end{equation}

\medskip The approach by Zhang and Zheng based on Young measures
follows two distinct routes.  Either one can use a viscous
regularization of the system \eqref{eq:wavesys} by adding the terms
$\epsilon R_{xx}$ and $\epsilon S_{xx}$ to the first and the second
equation, respectively, and subsequently analyze in detail the
behavior of the solution as $\epsilon\to0$, see \cite{ZZ:rarefactive}.
Alternatively \cite{ZZ:2001a,ZZ:2003:weakvariational,ZZ:2004}, one can
replace the quadratic growth on the right-hand side of equation
\eqref{eq:wavesys} by a linear growth for large values of $R^2$ and
$S^2$. More specifically, introduce the function
$$
Q_\epsilon(P)=\begin{cases}
\frac2\epsilon(\abs{P}-\frac{1}{2\epsilon}) & \text{for
  $\abs{P}\ge\frac1\epsilon$},\\  
 P^2 &\text{for $\abs{P}\le\frac1\epsilon$},
\end{cases}
$$
and replace the terms $R^2$ and $S^2$ by $Q_\epsilon(R)$ and
$Q_\epsilon(S)$, respectively, in the first and the second equation.
Again the behavior of the solution has to be analyzed carefully as
$\epsilon\to0$.

Our approach is based on Young measures for a semi-discrete finite
difference upwind scheme.  More precisely, introduce a positive
discretization parameter $\Dx$, and approximate $R(j\Dx,t)$ and
$S(j\Dx,t)$ by $R_j(t)$ and $S_j(t)$, respectively, that is,
$R(j\Dx,t)\approx R_j(t)$ and $S(j\Dx,t)\approx S_j(t)$, $j\in\Z$.
Observe that we keep the time variable continuous.  The dynamics of
$(R_j(t), S_j(t))$ is governed by the system of ordinary differential
equations
\begin{align*}
  R_j'(t) - \cj{+}(t)\Dp R_j(t) &=
  \tc_j(t)\left(R_j^2(t)-S_j^2(t)\right),\\
  S_j'(t) + \cj{-}(t)\Dm S_j(t) &=
  -\tc_j(t)\left(R_j^2(t)-S_j^2(t)\right),
\end{align*}
where $D_{\pm}K_j=\pm(K_{j\pm 1}-K_j)/\Dx$.  Furthermore, the
functions $c_{j\pm1/2}$ and $\tc_j$ are defined as functions of $R_j$
and $S_j$, cf.\ Section~\ref{sec:diffscheme}.  The system is augmented
by appropriate initial-data. We recover the function $u_{\Dx}$ by the
formula
\begin{equation*}
  \int_0^{u_{\Dx}(x,t)} 2c(u)\,du=\int^x (R_{\Dx}(\tilde
  x,t)-S_{\Dx}(\tilde x,t))\,d\tilde x, 
\end{equation*}
where $R_{\Dx}$ equals $R_j(t)$ on $[(j-\frac12)\Dx, (j+\frac12)\Dx)$,
and similarly for $S_{\Dx}$. We first show that the system possesses
solutions that are local in time, and a subsequent a priori estimate
turns the local solution into a global one.  Once the existence of
solutions of the ordinary differential equations has been established,
we follow the approach of Zhang and Zheng closely.

Formally, a smooth solution of \eqref{eq:wavesysCONS} will satisfy the identity
\begin{equation} \label{eq:wavesysCONSf}
(f(R)+f(S))_t-(c(u)(f(R)-f(S)))_x =2\tc  H(R,S),
\end{equation}
where
\begin{equation} \label{eq:wavesysCONSfAA}
H(R,S)=\frac12 (R^2-S^2)(f'(R)-f'(S)) -(f(R)-f(S))(R-S), 
\end{equation}
for any smooth function $f$. The corresponding discrete relation,
Lemma~\ref{lem:Hlem}, is rather more complicated. However, based on
this, one shows that the difference scheme keeps the $L^2$ norm of
$\{R_j,S_j\}$ from increasing, cf.~Corollary~\ref{cor:l2bound}; a
similar result holds in the continuous case as well, cf.~\cite[Lemma
1]{ZZ:rarefactive}.  Intrinsic to the equation is a blow-up property
that is not fully understood. Indeed it is known, see
\cite{GlasseyHunterZheng}, that there exist examples with $R$ and $S$
of opposite sign initially, where the solution becomes unbounded.
However, if the initial data both are negative initially, one can show
that the solution remains regular, see, e.g.,
\cite[Thm.~2]{ZZ:rarefactive}.  Henceforth we will make the assumption
here that $R$ and $S$ are nonpositive initially. As in the continuous
case, \cite[Lemma 4]{ZZ:rarefactive}, one can show also in the
discrete case, Lemma \ref{lem:invariant}, that the equation enjoys
invariant domains: If $(R_{\Dx},S_{\Dx})$ both are nonpositive
initially, then they will remain so. Furthermore, if
$(R_{\Dx},S_{\Dx})$ in addition are bounded from below initially, they
will remain so, with the same lower bound.
From this it follows that $L^p$ norms do not increase, cf.~\cite[Lemma
5]{ZZ:rarefactive} and Lemma \ref{lem:lpbounds}.  Using this one can
show in the discrete case, cf. Lemma \ref{lem:udxconverg}, using the
Arzelà--Ascoli theorem, that there exists a function $u$ such that
$$
\udx \to u \quad\text{uniformly on compacts in $\R\times [0,T]$.}
$$
The remaining part of the analysis is to show that the limit indeed
satisfies the equation.  From a priori $L^p$ bounds we infer that
$R_{\Dx}\weakstar \wlim{R}$ and $S_{\Dx}\weakstar \wlim{S}$ in
$L^\infty([0,T];L^2(\R))$ (recall that $R_{\Dx}$ equals $R_j(t)$ on
$[(j-\frac12)\Dx, (j+\frac12)\Dx)$, and similarly for $S_{\Dx}$), and
$(R_{\Dx}-S_{\Dx})^2\weakto \wlim{(R_{\Dx}-S_{\Dx})^2}$ in
$\Lenloc(\Pt)$.  Using the div-curl lemma, Lemma \ref{lem:divcurl},
and Murat's lemma, Lemma \ref{Murat}, we show that
$R_{\Dx}S_{\Dx}\weakto \wlim{R}\, \wlim{S}$ in $\Lenloc(\Pt)$,
cf.~Lemma \ref{lem:weakconvRS}.  Thus we have established that
\begin{equation*}
\left(\wlim{R}-\wlim{S}\right)_t -
\left(c(u)(\wlim{R}+\wlim{S})\right)_x = 0
\end{equation*}
holds weakly, cf.\ \eqref{eq:weakRS}.  
By direct analysis of the scheme we infer that
\begin{equation*}
  c(u)_x = 2\tc(u)(\wlim{R}-\wlim{S})\quad\text{weakly,}
\end{equation*}
cf.\ \eqref{eq:dclim}.  Using the weak identity
$(c(u)u_t)_x=(c(u)u_x)_t$ we infer that
$u_t=\frac12(\wlim{R}+\wlim{S})$ holds weakly.  To complete the
argument, we derive an evolution equation for $\wlim{R^2}+\wlim{S^2}$,
see Lemma \ref{lem:frenormlemma} (cf.\ \cite[Lemma 11]{ZZ:rarefactive})
to conclude that $u_{tt}-c(u)(c(u)u_x)_x=0$ weakly, and indeed that $u$ is a
weak solution.  Our main result can be stated as follows (cf.\ Theorem
\ref{thm:l3negconv}): If $u_0$ and $v_0$ are such that $R(\dott,0)$
and $S(\dott,0)$ are nonpositive, and in $L^3(\R)\cap L^1(\R)$, then
the semi-discrete difference scheme produces a sequence that converges
to a solution of \eqref{eq:waveeq} in the sense of Definition
\ref{def:weaksol}.
 
In Appendix~\ref{sec:higher} we show a higher integrability result,
see Lemma~\ref{lem:1alphabnd}, which reads, here stated in the
continuous case (cf.~\cite[Lemma 5]{ZZ:rarefactive}), as follows: \\ If
$(R(\dott,0),S(\dott,0))\in L^1(\R)\cap L^2(\R)$, then 
 $u_x\in L^{p}_{\mathrm{loc}}(\R\times [0,T],c'(u)dx)$ for $p\in
[2,3)$. The other results in this paper are independent of this, and
the significance of the appendix is that it is suspected that such a
regularity property could play a role in a uniqueness result.

\section{The difference scheme}
\label{sec:diffscheme}
Our first aim is to construct an approximate solution of
\eqref{eq:wavesys} based on a finite difference approximation of the
spatial derivative. Rather than work on the full system of three equations,
we derive approximate relations for the functions $c(u)$ and $\tc(u)$
in terms of $R$ and $S$, thereby reducing the system to two equations.
The temporal variable will not be discretized, and thus we will
consider systems of ordinary differential equations indexed by the
spatial lattice and depending on the lattice spacing.  Subsequently we
will show that as the lattice spacing decreases to zero, the system
converges to the solution of \eqref{eq:wavesys}. 

To avoid complicating the convergence analysis, we have chosen to
restrict our attention to a semi-discrete difference scheme.  To turn
the difference scheme into a fully discrete one, we can rely on a
variety of different time-discretization techniques, see Section
\ref{sec:numer} for one example.

We shall use \eqref{eq:wavesys} as a starting point for a difference
scheme. For $j\in\Z$, define $x_j=j\Dx$ and $x_{j+1/2}=x_j+\frac12\Dx$
where $\Dx>0$ is the lattice spacing. Let $I_j$ denote the interval
$[x_{j-1/2},x_{j+1/2})$.  

Given any function $K\colon \R\to\R$, we let the value of $K$ at the
point $x_j$ be denoted by $K_j$, that is, $K_j= K(x_j)$.  

On the other hand, given any sequence $\{K_j\}_{j\in\Z}$, we can
consider it as the sampling at lattice points $\Dx\, \Z$ of the
function $K$ defined by
\begin{equation}
K(x)=\sum_{j\in\Z} K_j \cha{I_j}(x). \label{eq:extend}
\end{equation}
Here $\cha{I}$ denotes the characteristic function of the set $I$.
Clearly, if values $\{K_j\}$ are computed from some difference scheme,
we consider the function \eqref{eq:extend} as the approximation of the
true solution.

It is easy to prove  the inequalities
$$
\norm{K}_{\el\infty}\le \frac{1}{\sqrt{\Dx}} \norm{K}_{\el2},\quad
\norm{K}_{\el2}\le \frac{1}{\sqrt{\Dx}}\norm{K}_{\el 1}.
$$
 Introduce forward and backward differencing by
\begin{equation*}
\Dpm K_j = \pm\frac{1}{\Dx}(K_{j\pm 1} - K_j)
\end{equation*}
for  any sequence $\seq{K_j}$ of real numbers.
Let $(R,S)=\{(R_j,S_j)\}_{j\in\Z}$
satisfy the (infinite) system of ordinary differential equations
\begin{align}
  R_j'(t) - \cj{+}(t)\Dp R_j(t) &=
  \tc_j(t)\left(R_j^2(t)-S_j^2(t)\right),\label{eq:Reqn}\\
  S_j'(t) + \cj{-}(t)\Dm S_j(t) &=
  -\tc_j(t)\left(R_j^2(t)-S_j^2(t)\right),\label{eq:Seqn}
\end{align}
for $j\in \Z$. 
The functions $\cj{\pm}$ and $\tc_j$ are specified as follows. First set
\begin{equation}\label{eq:Fdef}
  F(u)=\int_0^u 2c(v)\,dv.
\end{equation}
Since $c(u)>0$, we have $F'(u)>0$, and $F$ is therefore one-to-one. We
start by defining $\Fjmh$ by 
\begin{equation}
  \label{eq:Fmjhdef}
  \left.
    \begin{aligned}
      \lim_{j\to-\infty}\Fjmh&=0,\\
      \Dp\Fjmh &= R_j-S_j
      \end{aligned}\right\}
      \quad\text{or}\quad \Fjph =\Dx\sum_{i=-\infty}^j (R_i-S_i).
\end{equation}
Then we can define $\ujph$ by
\begin{equation}
  \label{eq:ujphdef}
  \ujph=(F^{-1})\left(\Fjph\right),\quad j\in\Z.
\end{equation}
Note that this implies
\begin{equation*}
  R_j-S_j=\Dp F(\ujmh)=2c\left(\bar{u}_j^+\right)\Dp\ujmh,
\end{equation*}
for some value $\bar{u}_j^+$ between $\uj{-}$ and $\uj{+}$.
Therefore
\begin{equation}
  \label{eq:ujdef2}
  \Dp \uj{-} = \frac{R_j-S_j}{2c\left(\bar{u}_j^+\right)}.
\end{equation}
Set 
\begin{equation}\label{eq:cjdef}
\cj{-}=c(\uj{-}),
\end{equation}
and note that for some $u_j^+$ between
$\uj{-}$ and $\uj{+}$ we have
\begin{equation}\label{eq:dc}
  \Dp\cj{-} = c'\left(u_j^+\right)\Dp\uj{-} =
  \frac{c'\left(u_j^+\right)}{2c\left(\bar{u}_j^+\right)} \left(R_j -
    S_j\right).
\end{equation}
So if we define
\begin{equation}\label{eq:tcdef}
\tc_j = \frac{c'\left(u_j^+\right)}{4c\left(\bar{u}_j^+\right)},
\end{equation}
we have that
\begin{equation}\label{eq:dccrucial}
\Dp \cj{-} = 2\tc_j \left(R_j - S_j\right).
\end{equation}
Thus we have defined the functions $\cj{\pm}=\cj{\pm}(R,S)$ and
$\tc_j= \tc_j(R,S)$.

We will work with $u_0\in H^1(\R)$ and $v_0\in L^2(\R)$.  In this case we define
\begin{equation}\label{eq:initapprox}
u_{0,j}=u_0(j\Dx), \quad u'_{0,j}=\frac1\Dx\int_{I_j}u'_0(x)\, dx, 
\quad v_{0,j}=\frac1\Dx\int_{I_j}v_0(x)\, dx.
\end{equation}
The initial values for \eqref{eq:Reqn} and \eqref{eq:Seqn} are
\begin{equation}
  R_j(0) = v_{0,j} + c\left(u_{j,0}\right) u_{0,j}',\quad\text{and}\quad
  S_j(0)= v_{0,j} - c\left(u_{j,0}\right) u_{0,j}', \label{RS_init}
\end{equation}
for $j\in\Z$.  Extend the initial data $\{(R_j(0),S_j(0))\}_{j\in\Z}$
by, cf.~\eqref{eq:extend}, 
\begin{align}\label{eq:start}
R_{0,\Dx}(x)&=R_{\Dx}(x,0)=\sum_j R_j(0)  \cha{I_j}(x), \notag \\
 S_{0,\Dx}(x)&=S_{\Dx}(x,0)=\sum_j S_j(0)  \cha{I_j}(x).
\end{align}
At this point it is convenient to state the following general lemma.
\begin{lemma}
  \label{lem:constl2converg} Let $\test$ be a function in $L^2(\R)$,
  and set 
  $$
  \test_j=\frac{1}{\Dx}\int_{I_j}\test(x)\,dx,\quad \test_{\Dx}(x)=\sum_j
  \test_j\cha{I_j}(x).
  $$
  Then 
  $$
  \norm{\test-\test_{\Dx}}_{L^2(\R)}\to 0
  $$
  as $\Dx\to 0$.
\end{lemma}
\begin{proof}
For general functions $\phi,\psi$ in $L^2(\R)$ we have
 \begin{align}
    \int_\R \left(\psi_{\Dx}(x)-\test_{\Dx}(x)\right)^2\,dx &=
    \sum_{j}\int_{I_j} \Bigl( \frac{1}{\Dx}\int_{I_j}
    (\psi(z)-\test(z))\,dz\Bigr)^2\,dx\notag \\
    &\le \sum_j \int_{I_j} \frac{1}{\Dx}\int_{I_j}
    \left(\psi(z)-\test(z)\right)^2\, dzdx \label{eq:dom0} \\
    &= \sum_j \int_{I_j} 
    \left(\psi(z)-\test(z)\right)^2 \,dz  \notag\\
    &= \int_\R \left(\psi(z)-\test(z)\right)^2\,dz.\notag
\end{align}
Thus
\begin{align}
\norm{\test-\test_{\Dx}}_2&\le \norm{\psi-\test}_2+\norm{\psi_{\Dx}-\test_{\Dx}}_2
+\norm{\psi-\psi_{\Dx}}_2\notag \\
&\le2 \norm{\psi-\test}_2+\norm{\psi-\psi_{\Dx}}_2
\end{align}
which shows that we, without loss of generality, can assume that $\test$ is a smooth function with compact support, say $\supp(\test)\subseteq[-M,M]$.  We find, similarly to the derivation of \eqref{eq:dom0}, that
  \begin{align}
    \int_\R \left(\test(x)-\test_{\Dx}(x)\right)^2\,dx &=
    \sum_{j}\int_{I_j} \Bigl( \frac{1}{\Dx}\int_{I_j}
    (\test(x)-\test(z))\,dz\Bigr)^2\,dx\notag \\
    &\le \sum_j \int_{I_j} \frac{1}{\Dx}\int_{I_j}
    \left(\test(x)-\test(z)\right)^2\, dzdx \label{eq:dom} \\
    &\le 2\int_{-M-1}^{M+1}\Biggl(\frac{1}{2\Dx}\int_{-\Dx}^{\Dx}
    \left(\test(x)-\test(x-y)\right)^2 \,dy \Biggr)\,dx. \notag
\end{align}
Since $\test$ is uniformly continuous, we can find $\delta>0$ such that 
$\abs{y}\le\delta$ implies
\begin{equation}
\abs{\test(x)-\test(x-y)}^2\le \frac{\eps}{4(M+1)}, \quad x\in \R.
\end{equation}
By choosing $\Dx\le\delta$ we find that
\begin{equation}
 \int_\R \left(\test(x)-\test_{\Dx}(x)\right)^2\,dx\le \eps,
\end{equation}
which concludes the proof. 
\end{proof}
This lemma implies
\begin{equation}\label{eq:start-sterkkonv}
  \begin{split}
    \norm{R_0-R_{0,\Dx}}_{L^2(\R)}\to 0 \quad \text{as $\Dx\to 0$,} \\
    \norm{S_0-S_{0,\Dx}}_{L^2(\R)}\to 0 \quad \text{as $\Dx\to 0$.}
  \end{split}
\end{equation}

\begin{hypothesis} \label{hyp:main}
  Consider $u_0\in W^{1,3+q}(\R)\cap W^{1,1}(\R)$ and $v_0\in
  L^{3+q}(\R)\cap L^1(\R)$ for some $q>0$, and let $R_0$ and $S_0$ be
  defined by \eqref{eq:R0S0def}.  Then we assume that $R_0\le 0$ and $S_0\le
  0$ almost everywhere.
\end{hypothesis}
This assumption implies that also $R_j(0)$ and $S_j(0)$ are
nonpositive for all $j$.  Furthermore, by interpolation, we have that
$u_0\in W^{1,p}(\R)$ and $v_0\in L^p(\R)$ for any $p\in[1,3+q]$.
 

\begin{lemma} \label{lem:global}
Assume Hypothesis \ref{hyp:main}. Then
  the system \eqref{eq:Reqn}--\eqref{eq:Seqn} of ordinary differential equations  with initial data
  \eqref{RS_init} has a unique $C^1$
  solution $\seq{R_j(t)}_{j\in\Z}$ and $\seq{S_j(t)}_{j\in\Z}$ for all
  $t>0$.
\end{lemma}
\begin{proof}
  We use the notation $R(t)=\seq{R_j(t)}_{j\in\Z}$ and
  $S(t)=\seq{S_j(t)}_{j\in\Z}$ 
  and write \eqref{eq:Reqn} and \eqref{eq:Seqn}, as 
  \begin{align*}
    R_j'(t)&=\Psi_j^R(R,S),\\
    S_j'(t)&=\Psi_j^S(R,S).
  \end{align*}
  Viewing this as an ordinary differential equation in
  $\ell^1(\Z)\times\ell^1(\Z)$, where $\ell^1(\Z)$ denotes the set of
  absolutely summable sequences with norm 
  $$
  \norm{v}_{\el1}=\Dx\sum_j\abs{v_j},
  $$
  it will have a unique  differentiable solution
  $(R(t),S(t))$ if $\Psi(R,S)=\seq{(\Psi^R_j,\Psi^S_j)}_{j\in\Z}$ is
  locally Lipschitz continuous. This solution will be defined for $t$
  in some interval $[0,t^*)$ where $t^*$ is a ``blow-up'' time, i.e., 
  $$
  \lim_{t\uparrow t^*} \big(\norm{R(t)}_{\el1}+\norm{S(t)}_{\el1}\big)=\infty.
  $$
  If one can show that $\norm{R(t)}_{\el1}+\norm{S(t)}_{\el1}<\infty$ for all
  $t>0$, then a continuously differentiable solution exists for all time. 

  Now we claim that
  \begin{equation}
    \label{eq:Psilip}
    \norm{\Psi(R,S)-\Psi(\hat{R},\hat{S})}_{\el1}\le
    C\left(\norm{R-\hat{R}}_{\el1}+\norm{S-\hat{S}}_{\el1}\right),
  \end{equation}
  where $C$ is a constant that depends on $\norm{(R,S)}_{\el1}$,
  $\|(\hat{R},\hat{S})\|_{\el1}$ and $\Dx$.  We shall show this for
  $\Psi^R$, the arguments for $\Psi^S$ are identical. 

  To show Lipschitz continuity we start by recalling (cf.~\eqref{eq:Fmjhdef})
  $$
   F_{j+1/2}=F_{j+1/2}(R,S)= \Dx\sum_{i=-\infty}^j (R_i-S_i).
  $$
  Then 
  \begin{align*}
    &\abs{F_{j+1/2}\left(R,S\right)-F_{j+1/2}(\hat{R},\hat{S})}\\
    &\qquad\qquad \le
    \Dx\sum_{i=-\infty}^j
    \left(\abs{R_i-\hat{R}_i}+\abs{S_i-\hat{S}_i}\right)\le
    \norm{R-\hat{R}}_{\el1}+\norm{S-\hat{S}}_{\el1}, 
  \end{align*}
  and therefore (writing $\hat{F}_j=F_{j+1/2}(\hat{R},\hat{S})$)
  \begin{equation}
    \label{eq:Philinfty}
    \norm{F-\hat{F}}_{\el\infty}\le 
    \norm{R-\hat{R}}_{\el1}+\norm{S-\hat{S}}_{\el1}. 
  \end{equation}
  Next we find (cf.~\eqref{eq:Reqn}) using \eqref{eq:ujphdef}, \eqref{eq:cjdef}, and 
  \eqref{eq:dccrucial} that
\begin{align} \label{eq:tildeLip}
\tc_j (R_j^2-S_j^2)&=\frac12 \Dp \cj{-}(R_j+S_j) \notag \\
&= \frac12(\cj{+}-\cj{-})(R_j+S_j) \notag \\
&= \frac12(c(F_{j+1/2})-c(F_{j-1/2}))(R_j+S_j),
\end{align}
abbreviating $c(F_{j\pm1/2})=c((F^{-1})(F_{j\pm1/2}))$.  Thus
\begin{align*}
\tc_j (R_j^2-S_j^2)- \hat{\tc}_j (\hat{R}_j^2-\hat{S}_j^2)&= 
 \frac12\Big(\big(c(F_{j+1/2})-c(\hat{F}_{j+1/2}) \big) (R_j+S_j) \\
 &\qquad\qquad -\big(c(F_{j-1/2})-c(\hat{F}_{j-1/2}) \big)(R_j+S_j) \\
&\quad+(c(\hat{F}_{j+1/2})-c(\hat{F}_{j-1/2}))(R_j- \hat{R}_j+S_j-\hat{S}_j) \Big).
\end{align*}
  Now 
  \begin{align*}
    \abs{\Psi^R_j(R,S)-\Psi^R_j(\hat{R},\hat{S})}&\le
    \abs{c\left(F_{j+1/2}\right)-c(\hat{F}_{j+1/2})}
    \abs{\Dp  R_j}\\
    &\qquad +
    c(\hat{F}_{j+1/2})\abs{\Dp\left(\hat{R}_j-R_j\right)}\\
 &\qquad +\frac12\abs{c(F_{j+1/2})-c(\hat{F}_{j+1/2})} (\abs{R_j}+\abs{S_j}) \\
 &\qquad\qquad +\frac12\abs{c(F_{j-1/2})-c(\hat{F}_{j-1/2})} (\abs{R_j}+\abs{S_j})  \\
&\quad+\frac12\abs{c(\hat{F}_{j+1/2})-c(\hat{F}_{j-1/2})}
\big(\abs{R_j- \hat{R}_j}+\abs{S_j-\hat{S}_j}\big)  \\  
    &\le \frac{C}{\Dx}\abs{F_{j+1/2}-\hat{F}_{j+1/2}}\left(\abs{R_j}+\abs{R_{j+1}}\right)\\
    &\qquad + \frac{C}{\Dx}\left(\abs{R_j-\hat{R}_j}+
      \abs{R_{j+1}-\hat{R}_{j+1}}\right)\\
 &\qquad + \frac{C}{2\Dx}\abs{F_{j+1/2}-\hat{F}_{j+1/2}}   (\abs{R_j}+\abs{S_j})  \\
&\qquad + \frac{C}{2\Dx}\abs{F_{j-1/2}-\hat{F}_{j-1/2}}   (\abs{R_j}+\abs{S_j}) \\
&\qquad + C
    \left(\abs{R_j-\hat{R}_j}+\abs{S_j-\hat{S}_j}\right)
  \end{align*}
  since  $c$ is Lipschitz continuous functions of $F_{j+1/2}$. 
  Multiplying the above by $\Dx$ and summing over $j$, we see that
  $$
  \norm{\Psi^R(R,S)-\Psi^R(\hat{R},\hat{S})}_{\el1}\le
  C\left(\norm{R-\hat{R}}_{\el1}+\norm{S-\hat{S}}_{\el1}\right),
  $$
  where we have used \eqref{eq:Philinfty} to find 
  $$
  C=C\left(\Dx,\norm{R}_{\el1},\norm{S}_{\el1},
  \norm{\hat{R}}_{\el1},\norm{\hat{S}}_{\el1}\right).
  $$
  Therefore \eqref{eq:Psilip} holds, and
  we have established that $\seq{R_j(t)}_{j\in\Z}$ and
  $\seq{S_j(t)}_{j\in\Z}$ exist for $t<t^*$ (for any initial data). If
  the initial data are nonpositive and in $\el1$,
  Lemma~\ref{lem:lpbounds} concludes the proof.
\end{proof}
\begin{remark} \label{rem:global}
The existence of global solutions of the system \eqref{eq:Reqn}--\eqref{eq:Seqn} with initial data 
\eqref{eq:initapprox}, that is, the fact that $t^*=\infty$, follows only after the estimate in 
Lemma~\ref{lem:lpbounds}, i.e., the inequality \eqref{eq:lpbounds}. Thus the results up to 
Lemma~\ref{lem:lpbounds} are first valid 
for all times less than $t^*$, and only after Lemma~\ref{lem:lpbounds} we can infer that  $t^*=\infty$. To simplify the notation, we state all these result for all $t$.  

For Lemma \ref{lem:global} we only require that $u_0\in W^{1,1}(\R)$ and $v_0\in L^1(\R)$.
\end{remark}

\section{Convergence analysis}\label{sec:converg}
Now let $f$ be a sufficiently smooth function, and observe
that
\begin{equation*}
f\left(R_{j+1}\right)=f\left(R_j\right) +
f'\left(R_j\right)\left(R_{j+1}-R_j\right) +
\frac{1}{2}f''(r_j)\left(R_{j+1}-R_j\right)^2,
\end{equation*}
where $r_j$ is between $R_{j+1}$ and $R_j$. This can be rewritten
\begin{equation}
  \label{eq:Dfr}
  \Dp f\left(R_j\right) = f'\left(R_j\right) \Dp R_j + \frac{\Dx}{2} 
  f''\left(r_j\right) \left(\Dp R_j\right)^2.
\end{equation}
Furthermore, we have for any quantity $f_j$,
\begin{align}
  \Dp\left(\cj{-}f_j\right) &= \cj{+}\Dp f_j + f_j \Dp\cj{-}\notag\\
  &=\cj{+} \Dp f_j + 2 f_j \tc_j \left(R_j-S_j\right).\label{eq:DfrC}
\end{align}
 Similarly to \eqref{eq:Dfr} we have for a
  sufficiently smooth function $g$
  \begin{equation}
    \label{eq:Dfs}
    \Dm g\left(S_j\right) = g'\left(S_j\right) \Dm S_j - \frac{\Dx}{2}
    g''\left(s_j\right)\left(\Dm S_j\right)^2,
  \end{equation}
  where $s_j$ is between $S_j$ and $S_{j-1}$. We also have
  \begin{align}
    \Dm\left(\cj{+}g_j\right) &= \cj{-} \Dm g_j + g_j \Dp\cj{-}\notag\\
    &= \cj{-}\Dm g_j + 2 g_j \tc_j \left(R_j - S_j\right).  \label{eq:DfsC}
  \end{align}
Next define (cf.~\eqref{eq:wavesysCONSfAA})
\begin{equation*}
H(R,S)=\frac{1}{2}\left(f'\left(R\right) -
    f'\left(S\right)\right)\left(R^2 - S^2\right) -
  \left(f\left(R\right)-f\left(S\right)\right)
  \left(R-S\right).
\end{equation*}
We shall use the following lemma repeatedly.
\begin{lemma}\label{lem:Hlem}  
Let $f,g\in C^2(\R)$.  Consider sequences $\{R_j\}_{j\in \Z}$, $\{S_j\}_{j\in \Z}$ satisfying \eqref{eq:Reqn}--\eqref{eq:Seqn}. 
For $t>0$, there holds
  \begin{align}
    \label{eq:fRcons11}
     \begin{aligned}
    \frac{d}{dt} f(R_j) -& \Dp\left(\cj{-} f(R_j)\right)
    + \frac{\Dx}{2} \cj{+} f''\left(r_j\right) \left(\Dp R_j\right)^2 \\
  &\qquad  = 2\tc_j \left( \frac{1}{2}f'(R_j)\left(R_j^2 - S_j^2\right) - f(R_j)
      \left(R_j -S_j\right)\right), 
      \end{aligned}\\
   \label{eq:fScons11}
     \begin{aligned}
    \frac{d}{dt} g(S_j) +& \Dm\left(\cj{+}g(S_j)\right) + \frac{\Dx}{2}\cj{-}
    g''\left(s_j\right) \left(\Dm S_j\right)^2 \\ 
    &\qquad = 
    -2\tc_j\left(\frac{1}{2}g'(S_j) \left(R_j^2 - S_j^2\right) 
      - g(S_j)\left(R_j - S_j\right)\right).
   \end{aligned}    
  \end{align}
In particular,
\begin{equation}
    \label{eq:sumcons}
    \begin{aligned}
      \frac{d}{dt} \bigl(f\left(R_j\right)+&f\left(S_j\right)\bigr)
      -\Dp\left(\cj{-} f\left(R_j\right)\right)
      +\Dm\left(\cj{+} f\left(S_j\right)\right)\\
      &\quad + \frac{\Dx}{2}\left(\cj{+}f''\left(r_j\right) \left(\Dp
          R_j\right)^2 + \cj{-} f''\left(s_j\right) \left(\Dm
          S_j\right)^2\right) \\
      &\qquad\qquad=2\tc_j H\left(R_j,S_j\right),
    \end{aligned}
\end{equation}
where $r_j$ is between $R_j$ and $R_{j+1}$, and $s_j$ between $S_j$
and $S_{j-1}$. 
\end{lemma}
\begin{proof}
  Multiplying \eqref{eq:Reqn} by $f'_j=f'(R_j)$, using \eqref{eq:Dfr} and \eqref{eq:DfrC}, we
  find that
  \begin{multline}
    \label{eq:fRcons}
    \frac{d}{dt} f_j - \Dp\left(\cj{-} f_j\right)
    + \frac{\Dx}{2} \cj{+} f''\left(r_j\right) \left(\Dp R_j\right)^2 
    \\= 2\tc_j \left( \frac{1}{2}f_j'\left(R_j^2 - S_j^2\right) - f_j
      \left(R_j -S_j\right)\right),
  \end{multline}
  where $f_j=f(R_j)$.
  Similarly, multiplying \eqref{eq:Seqn} with $g'_j=g'(S_j)$ for some function $g\in C^2(\R)$, using
  \eqref{eq:Dfs} and \eqref{eq:DfsC}, we find that
  \begin{multline}
    \label{eq:fScons}
    \frac{d}{dt} g_j + \Dm\left(\cj{+}g_j\right) + \frac{\Dx}{2}\cj{-}
    g''\left(s_j\right) \left(\Dm S_j\right)^2 \\ = 
    -2\tc_j\left(\frac{1}{2}g_j' \left(R_j^2 - S_j^2\right) 
      - g_j\left(R_j - S_j\right)\right),
  \end{multline}
  where $g_j=g(S_j)$.  Choosing $f=g$ and adding \eqref{eq:fRcons} and
  \eqref{eq:fScons} we conclude that \eqref{eq:sumcons} holds.
\end{proof}
This lemma has several useful consequences, the first of which is the following result.
\begin{corollary}\label{cor:l2bound}
We have that
\begin{multline}
      \Dx\sum_j \left(R_j^2 + S_j^2\right) (t)  +\int_0^T \Dx \sum_j
      \Dx \left(\cj{+}\left(\Dp R_j\right)^2 + \cj{-}\left(\Dm
      S_j\right)^2\right)\, dt\\ 
      \le \Dx\sum_j \left(R_j^2 + S_j^2\right) (0).
      \label{eq:l2bound}
\end{multline}
In particular, we have
\begin{equation} \label{eq:l2bound0}
        \Dx\sum_j \left(R_j^2 + S_j^2\right) (t) 
        \le \Dx\sum_j \left(R_j^2 + S_j^2\right) (0).
\end{equation}
\end{corollary}

\begin{proof}
  Appy Lemma~\ref{lem:Hlem} with $f(K)=K^2$. In this case we observe
  that $H(R,S)=0$, and  $f''=2$.  Therefore,
  Lemma~\ref{lem:Hlem} yields
  \begin{align*}
    \frac{d}{dt}\left(R^2_j+S^2_j\right) &- \Dp\left(c_{j-1/2}R_j^2\right) +
    \Dm\left(c_{j+1/2}S_j^2\right)\\
    &+ \Dx  \left(\cj{+}\left(\Dp
        R_j\right)^2 
      + \cj{-}\left(\Dm S_j\right)^2\right)\le 0.
  \end{align*}
  Multiplying with $\Dx$, summing over $j$, and integrating in $t$
  finishes the proof of the corollary.
\end{proof}
The variational wave equation enjoys certain invariance properties in
the $(R,S)$ variables.  Indeed, if both are nonpositive initially,
they will remain so for all time. Furthermore, if in addition the
initial data are bounded below by a (negative) constant, then the same
constant bounds the solution for all time.  See, e.g., \cite[Thm.
3.1.6]{ZZ:2005}.  The approximate solution has the same properties,
which is the result of the following lemma.
\begin{lemma}
  \label{lem:invariant}
  The following statements hold: \\
 (i)     If $R_j(0)\le 0$ and $S_j(0)\le 0$ for all $j$, then $R_j(t)\le 0$
      and $S_j(t)\le 0$ for all $t\ge 0$ and all $j$. \\      
  (ii)    If $-M\le R_j(0)\le 0$ and $-M\le S_j(0)\le 0$ for some positive
      number $M$ and for all $j$, then 
      $-M\le R_j(t)\le 0$ and $-M\le S_j(t)\le 0$ for all $j$ and
      $t\ge 0$.
\end{lemma}
\begin{proof}
  To prove the first statement (\textit{i}) choose
  $f(K)=\left(\max\seq{0,K}\right)^2$ in \eqref{eq:sumcons}; with this choice
  \begin{equation*}
  H(R,S)=
  \begin{cases}
    0 \quad &\text{if $RS\ge 0$,}\\
    RS(R-S) &\text{if $S<0$ and $R>0$,}\\
    RS(S-R) &\text{if $S>0$ and $R<0$.}
  \end{cases}
  \end{equation*}
  Hence $H(R,S)\le 0$, furthermore $f''(K)\ge 0$, and thus using
  \eqref{eq:sumcons} we find that
  \begin{equation*}
  \sum_j \big(\left(\max\seq{0,R_j(t)}\right)^2 +
  \left(\max\seq{0,S_j(t)}\right)^2\big) \le 0,
  \end{equation*}
  since $R_j(0)\le 0$ and $S_j(0)\le 0$ for all $j$. Thus the first
  statement (\textit{i})  of the lemma holds. 
  
  To prove the second statement (\textit{ii}) choose 
  $f(K)=\left(\min\seq{K+M,0}\right)^2$. Then we find that
  \begin{equation*}
  H(R,S) = 
  \begin{cases}
    0 \quad &\text{if $R\ge -M$ and $S\ge -M$,}\\
    -2M(R-S)^2 &\text{if $R<-M$ and $S<-M$,}\\
    (R+M)(R-S)(S-M) &\text{if $R<-M\le S$,}\\
    (S+M)(R-S)(M-R) &\text{if $S<-M\le R$,}
  \end{cases}
  \end{equation*}
  which implies  that 
  \begin{equation}\label{eq:HRS}
  H(R,S)\bigm|_{ \seq{(R,S)\mid R<M,S<M}} \le 0.
  \end{equation}
  Furthermore $f''(K)\ge 0$.  Thus, if $0\ge R_j(0)\ge -M$ and $0\ge
  S_j(0)\ge -M$, we observe from the first statement (\textit{i}) that
  $R_j$ and $S_j$ remain negative. This implies, using \eqref{eq:HRS}, that 
  $H(R_j,S_j)(t)\le 0$. Hence it follows as before, using equation
  \eqref{eq:sumcons}, that
  \begin{equation*}
    \sum_j \big(\left(\min\seq{0,R_j(t)+M}\right)^2 +
    \left(\min\seq{0,S_j(t)+M}\right)^2\big) \le 0.
  \end{equation*}
  Thus the second statement (\textit{ii}) of the lemma follows.
\end{proof}
In case Hypothesis~\ref{hyp:main} holds, we have the integrability estimate.
\begin{lemma}
\label{lem:lpbounds} 
If $R_j(0)\le 0$ and $S_j(0)\le 0$ for all $j$, then
\begin{equation}\label{eq:lpbounds}
  \Dx\sum_j\big( \abs{R_j(t)}^p + \abs{S_j(t)}^p \big)\le 
  \Dx \sum_j \big(\abs{R_j(0)}^p + \abs{S_j(0)}^p\big),
\end{equation}
for any  $p\ge 1$. In addition, if Hypothesis~\ref{hyp:main} holds,
then for any $p\in (2,3+q)$
\begin{equation}
  \label{eq:highint2}
  \int_0^T \Dx\sum_j c'\left(u_j^+\right)
  \abs{\Dp \uj{-}}^{p+1}
  \,dt \le C_{p,T},
\end{equation}
where $C_{p,T}$ is a constant depending on $p$ and $T$ (but not on $\Dx$).  
\end{lemma}
\begin{remark} This lemma finishes the proof of
  Lemma~\ref{lem:global}, namely the fact that $t^*=\infty$,
  cf.~Remark~\ref{rem:global}.  
\end{remark}
\begin{proof}[Proof of {Lemma~\ref{lem:lpbounds}}]
  Choose $f(K)=\abs{K}^p$, and observe that
  \begin{equation*}
    \text{$f(0)=0$, $f(K)=f(-K)$ and $f''(K)\ge 0$}.
  \end{equation*}
  Now it is easy to see that
  \begin{equation*}
    H(R,S)=H(S,R)\quad\text{and}\quad H(-S,-R)=-H(R,S).
  \end{equation*}
  Furthermore
  \begin{equation*}
    H(R,R)=H(R,-R)=0.
  \end{equation*}
  We also find that
  \begin{equation*}
  \nabla H(R,S)\cdot(1,1) = \frac{1}{2}\left(f''(R)-f''(S)\right)
  \left(R^2-S^2\right)\ge 0,
  \end{equation*}
  since $f''$ is an even non-negative function. From this it follows that
  \begin{equation*}
    H(R,S)\bigm|_{\seq{(R,S)\mid R+S\le 0}}\le 0.
  \end{equation*}
  Hence, since $S_0\le 0$ and $R_0\le 0$, by Lemma~\ref{lem:invariant}
  also $R_j(t)$ and $S_j(t)$ are nonpositive for $t>0$, and thus
  \begin{equation*}
    \Dx\sum_j \big(f(R_j(t))+f(S_j(t)) \big) \le \Dx\sum_j
    \big(f(R_j(0)) + f(S_j(0))\big). 
  \end{equation*}
  
  \medskip
  For  the proof of \eqref{eq:highint2}, we fix $p\in (2,3+q)$,
  remember that $R_j\le 0$ and $S_j\le 0$, and
  calculate
  \begin{align*}
    H\left(R_j,S_j\right)&=\frac{1}{2}\Bigl[
      p\left(\sgn{R_j}\abs{R_j}^{p-1}-\sgn{S_j}\abs{S_j}^{p-1}\right)
      \left(R_j^2-S_j^2\right)\\
      &\qquad - 2\left(\abs{R_j}^p-\abs{S_j}^p\right)
      (R_j-S_j)\Bigr]\\
      &=\frac{1}{2}\Bigl[
      -p\left(\abs{R_j}^{p-1}-\abs{S_j}^{p-1}\right)
      \left(\abs{R_j}^2-\abs{S_j}^2\right) \\
      &\qquad +
      2\left(\abs{R_j}-\abs{S_j}\right)\left(\abs{R_j}^p-\abs{S_j}^p\right)
      \Bigr]\\
      &=\frac{1}{2}\Bigl[-p
      \left(\abs{R_j}-\abs{S_j}\right)^2
      \left(\abs{R_j}^{p-1}+\abs{S_j}^{p-1}\right)\\
      &\qquad
      \underbrace{-2(p-1)\abs{R_j}\abs{S_j}\left(\abs{R_j}-\abs{S_j}\right)
        \left(\abs{R_j}^{p-2}-\abs{S_j}^{p-2}\right)}_{b(R_j,S_j)}\\
      &\qquad -2\abs{R_j}\abs{S_j}\left(\abs{R_j}-\abs{S_j}\right)
      \left(\abs{R_j}^{p-2}-\abs{S_j}^{p-2}\right) \\
      &\qquad +
      2\left(\abs{R_j}-\abs{S_j}\right)\left(\abs{R_j}^p-\abs{S_j}^p\right) 
      \Bigr]\\
      &=\frac{1}{2}\Bigl[
      -p
      \left(\abs{R_j}-\abs{S_j}\right)^2
      \left(\abs{R_j}^{p-1}+\abs{S_j}^{p-1}\right)+b(R_j,S_j)\\
      &\qquad +2\left(\abs{R_j}-\abs{S_j}\right)
      \left(-\abs{R_j}\abs{S_j}
        \left(\abs{R_j}^{p-2}-\abs{S_j}^{p-2}\right)
        +\abs{R_j}^p-\abs{S_j}^p\right)  
      \Bigr]\\
      &=\frac{1}{2}\Bigl[
      -p
      \left(\abs{R_j}-\abs{S_j}\right)^2
      \left(\abs{R_j}^{p-1}+\abs{S_j}^{p-1}\right)+b(R_j,S_j)\\
      &\qquad +2\left(\abs{R_j}-\abs{S_j}\right)^2
      \left(\abs{R_j}^{p-1}+\abs{S_j}^{p-1}\right)
      \Bigr]\\
      &=\frac{1}{2}\Bigl[
      -(p-2)
      \left(\abs{R_j}-\abs{S_j}\right)^2
      \left(\abs{R_j}^{p-1}+\abs{S_j}^{p-1}\right)+b(R_j,S_j)\Bigr].
  \end{align*}
  It is easy to see that $b(R,S)\le 0$ for $p>2$, and we also have the
  inequality
  $$
  \abs{R}^{p-1}+\abs{S}^{p-1}\ge K_{p}
  \abs{\left(\abs{R}-\abs{S}\right)}^{p-1},
  $$
  for some positive constant $K_p$ depending on $p$. Hence
  $$
  H\left(R_j,S_j\right)\le
  \frac{-K_p(p-2)}{2}\abs{R_j-S_j}^{p+1}=-K_p(p-2)c\left(\bar{u}_j^+\right)
  \abs{\Dp \uj{-}}^{p+1}.
  $$
  By Hypothesis~\ref{hyp:main}, we find that
  $$
  \frac{K_p(p-2)}{4}\int_0^T \Dx \sum_j c'\left(u^+_j\right) 
  \abs{\Dp \uj{-}}^{p+1} \le \Dx\sum_j
  \left(\abs{R_j(0)}^{p}+\abs{S_j(0)}^p\right) \le C,
  $$
  from which \eqref{eq:highint2} follows.
\end{proof}
Extend the functions $(R_j,S_j)$ to the full line,
cf.~\eqref{eq:extend} and \eqref{eq:start}, by
\begin{equation}\label{eq:RSdxdef}
  R_{\Dx}(x,t)=\sum_j R_j(t) \cha{I_j}(x),\quad\text{and}\quad 
  S_{\Dx}(x,t)=\sum_j S_j(t) \cha{I_j}(x).
\end{equation}
Define $F_{\Dx}$ by
\begin{equation}
  \label{eq:Fdxdef}
  F_{\Dx}(x,t)=\int^x (R_{\Dx}(\tilde x,t)-S_{\Dx}(\tilde x,t))\,d\tilde x,
\end{equation}
and then $u_{\Dx}$ by 
\begin{equation}
  \label{eq:udxdef}
  \int_0^{u_{\Dx}(x,t)} 2c(u)\,du = F_{\Dx}(x,t).
\end{equation}
Note that 
\begin{equation*}
  \Dp F_{\Dx}\left(x_{j-1/2},t\right)=R_j-S_j=\Dp F\left(u_{j-1/2}(t)\right),
\end{equation*}
or
\begin{equation*}
  \int_{u_{j-1/2}}^{u_{j+1/2}} 2c(v)\,dv = 
  \int_{u_{\Dx}(x_{j-1/2},t)}^{u_{\Dx}(x_{j+1/2},t)}\!\! 2c(v)\,dv.
\end{equation*}
Now we have that $\lim_{x\to-\infty}\Rdx(x,t)=\lim_{x\to-\infty}
\Sdx(x,t)=0$, and therefore $\lim_{x\to -\infty}\udx(x,t)=0$. Hence we
must have $\udx(x_{j-1/2},t)=u_{j-1/2}(t)$ for all $j$. It is
convenient also to define the piecewise constant function
\begin{equation}
  \label{eq:ubardef}
  \bar{u}_{\Dx}(x,t)=\sum_j u_{j-1/2}(t)
  \cha{I_{j-1/2}} (x).
\end{equation}
Now we can show the (local) uniform convergence of $\udx$.
\begin{lemma}
  \label{lem:udxconverg}  Assume Hypothesis \ref{hyp:main}.
  Then there
  exists a function $u\in C(\Pt)$
  such that for any finite interval $[a,b]$ we have
  \begin{equation*}
    \lim_{\Dx\to 0}\udx(x,t)=u(x,t) \quad\text{uniformly for $(x,t)\in
    [a,b]\times[0,T]$.}  
\end{equation*}
\end{lemma}
\begin{proof}
From Hypothesis \ref{hyp:main} we infer that
\begin{equation}\label{eq:initbnd}
  \norm{R_{\Dx}(\dott,0)}_{L^p(\R)} +
  \norm{S_{\Dx}(\dott,0)}_{L^p(\R)} \le C,
\end{equation}
for both $p=1$ and $p=2$  for some constant $C$ that is independent of $\Dx$.  
From Lemma \ref{lem:lpbounds} it follows that 
\begin{equation} \label{eq:time1}
  \norm{\Rdx(\dott,t)}_{L^1(\R)} + \norm{\Sdx(\dott,t)}_{L^1(\R)} \le
  \norm{\Rdx(\dott,0)}_{L^1(\R)} + \norm{\Sdx(\dott,0)}_{L^1(\R)} \le C
\end{equation}
where $C$ is the constant in \eqref{eq:initbnd}. This implies that
$F_{\Dx}$ is uniformly bounded, since
$$
\abs{F_{\Dx}(x,t)} = \abs{\int^x (\Rdx(y,t)-\Sdx(y,t))\,dy} \le C.
$$
Next, we observe that
$$
F_{\Dx}(x,t)-F_{\Dx}(0,t)=\int_0^x (\Rdx(y,t)-\Sdx(y,t))\,dy.
$$
Therefore, using Cauchy--Schwarz's inequality and \eqref{eq:time1}, we find 
\begin{align*}
  \norm{F_{\Dx}(\dott,t)-F_{\Dx}(0,t)}_{L^2(a,b)}^2 &=
  \int_a^b \biggl( \int_0^x (\Rdx(y,t)-\Sdx(y,t))\,dy\biggr)^2\,dx\\
  &\le \int_a^b 2\abs{x} \int_{\R} (\Rdx^2(y,t)+\Sdx^2(y,t))\dy \,dx\\
    &\le C^2 \left(a^2+b^2\right).
  \end{align*}
   Similarly, by using
  \eqref{eq:Fdxdef} we find that 
   \begin{align*} 
  \norm{\partial_x F_{\Dx}(\dott,t)}_{L^2(a,b)}^2 &\le \int_a^b \bigg(\Rdx(x,t)-\Sdx(x,t) \bigg)^2\, dx\\
&  \le 2 \int_a^b \bigg(\abs{\Rdx(x,t)}^2+\abs{\Sdx(x,t)}^2 \bigg)\, dx \\
&\le    2(a+b)C^2
  \end{align*}
 using \eqref{eq:l2bound0}.
  Thus there is a constant $C_1$, independent of $t$ and $\Dx$ (but depending on $a,b$), such that 
  \begin{equation*}
  \norm{F_{\Dx}(\dott,t)-F_{\Dx}(0,t)}_ {H^1(a,b)}\le C_1.
   \end{equation*}
  Morrey's inequality now implies that for $x$ and $y$ in $[a,b]$ we have that 
  \begin{equation}
    \label{eq:Fcontx}
    \abs{F_{\Dx}(x,t)-F_{\Dx}(y,t)}\le C_2 \abs{x-y}^{1/2},
  \end{equation}
  for some constant $C_2$ which is independent of $\Dx$ and $t$ (but depending on $a,b$). 

  Now note that using $f(R)=R$ in \eqref{eq:fRcons} and
  $g(S)=S$ in \eqref{eq:fScons} we find, cf.~\eqref{eq:wavesysCONS}, that
  \begin{align*}
    R_j' -
    \Dp\left(\cj{-}R_j\right)&=-{\tc_j}\left(R_j-S_j\right)^2,\\
    S_j' + \Dm\left(\cj{+}S_j\right)&= -{\tc_j}\left(R_j-S_j\right)^2.
  \end{align*}
  Since $\Dp F_j = R_j -S_j$,
  \begin{equation*}
  \frac{d}{dt} F_j = \Dx \sum_{i=-\infty}^{j-1} (R_i' - S_i') = \cj{-}\left(R_{j}
    + S_{j-1}\right).
  \end{equation*}
  Therefore we have that for $0\le s\le t\le T$
  \begin{align*}
    \norm{F_{\Dx}(\dott,t)-F_{\Dx}(\dott,s)}_{L^2(\R)} &\le
    \int_s^t \norm{\frac{d}{dt} F_{\Dx}(\dott,\tau)}_{L^2(\R)}\,d\tau
    \\
    &\le 
    C_3\int_s^t \bigg(\norm{S_{\Dx}(\dott,\tau)}_{L^2(\R)} +
    \norm{R_{\Dx}(\dott,\tau)}_{L^2(\R)}\bigg) \, d\tau\\
    &\le C_3 C \abs{t-s},
  \end{align*}
using \eqref{eq:l2bound0}. 

Since  $H^1(a,b)\subset\!\subset C(a,b) \subset L^2(a,b)$, we can use
 \cite[Lemma 8]{simon}
 to deduce that for $x$ and $y$ in $(a,b)$, we
  have that for any $\eta>0$, there is a finite $C_\eta>0$ such that 
  \begin{align*}
    \abs{F_{\Dx}(x,t)-F_{\Dx}(x,s)}&\le
    \eta \norm{F_{\Dx}(\dott,t)-F_{\Dx}(\dott,s)}_{H^1(a,b)}\\
    &\qquad
    + 
    C_{\eta}\norm{F_{\Dx}(\dott,t)-F_{\Dx}(\dott,s)}_{L^2(a,b)}\\
    &\le \eta 2 C_1 + C_{\eta} C C_3 \abs{t-s}.
  \end{align*}
  For any $\epsilon>0$ we
  choose $(x,t)$ and $(y,s)$ in $[a,b]\times [0,T]$ and $\eta>0$ such that 
  $$
  C_2\abs{x-y}^{1/2}\le \frac{\epsilon}{3}, \quad
  \eta 2 C_1 \le \frac{\epsilon}{3} \quad\text{and then}\quad
  C_\eta C C_3 \abs{t-s} \le \frac{\epsilon}{3}.
  $$
  With this choice
  \begin{align*}
    \abs{F_{\Dx}(x,t)-F_{\Dx}(y,s)} &\le 
    \abs{F_{\Dx}(x,t)-F_{\Dx}(y,t)} + \abs{F_{\Dx}(y,t)-F_{\Dx}(y,s)}
    \\
    &\le \epsilon.
  \end{align*}
  Hence, the sequence $\seq{F_{\Dx}}_{\Dx>0}$ is equicontinuous and
  uniformly bounded in
  $[a,b]\times[0,T]$, and by the Arzelà--Ascoli theorem there exists a
  convergent subsequence (which we do not relabel).

  By the definition \eqref{eq:udxdef} of $u_{\Dx}$ and the assumption on $c$, 
  cf.\ \eqref{c_est}, we find that
  \begin{equation*}
    \abs{F_{\Dx_j}(x,t)-F_{\Dx_k}(x,t)}\ge C_4
    \abs{u_{\Dx_j}(x,t)-u_{\Dx_k}(x,t)},
  \end{equation*}
  for some constant $C_4$ depending only on the function $c$.
  This shows that $\seq{u_{\Dx_j}}$ is Cauchy and thus uniformly
  convergent on compacts $[a,b]\times[0,T]$.
\end{proof}
\begin{remark}
  For this lemma to hold it is sufficent to assume that $R_0$ and
  $S_0$ (and therefore $\Rdx(0),
  \Sdx(0)$) are nonpositive,  and in $L^1\cap L^2$.
\end{remark}

Note that $F_{\Dx}$ is linear in the interval $I_j$ as $R_{\Dx}$
and $S_{\Dx}$ are constant there. By definition we have that 
$$
\frac{\partial \udx}{\partial F_{\Dx}} =\frac{1}{2c(\udx)}>0.
$$
This means that $\udx(x,t)$ is monotone in the interval $I_j$, and we
have  that $\udx(x_{j\pm 1/2},t)=u_{j\pm 1/2}(t)$. 
To simplify the subsequent calculations we introduce
$$
\tilde{u}_j = \theta_j u_{j-1/2}+(1-\theta_j)u_{j+1/2}, \quad
\theta_j\in [0,1],\quad j\in \Z,
$$
and  
$$
\tilde{u}_{\Dx}(x,t)=\sum_j \tilde{u}_j(t) \cha{I_j}(x).
$$
Then for any fixed $x$ and $t$,
\begin{equation}\label{eq:tildelim}
  \lim_{\Dx\to 0}\tilde{u}_{\Dx}(x,t)=u(x,t).
\end{equation}
This is so since if $x\in I_j$, there is a $y_j\in I_j$ such that
$\tilde{u}_{\Dx}(x,t)=\udx(y_j,t)$ by the monotonicity of $\udx$. Now
let $\seq{\Dx_k}$ and $\seq{\Dx_\ell}$ be two sequences tending to
zero. Fixing $x$, we can find sequences $\seq{y_k}$ and $\seq{y_\ell}$
such that $y_k\to x$ and $y_\ell\to x$ and
$$
\tilde{u}_{\Dx_k}(x,t)=u_{\Dx_k}\left(y_k,t\right)\quad\text{and}\quad
\tilde{u}_{\Dx_\ell}(x,t)=u_{\Dx_\ell}\left(y_\ell,,t\right).
$$
Hence
$$
\abs{\tilde{u}_{\Dx_k}(x,t)-\tilde{u}_{\Dx_l}(x,t)}\le 
\abs{u_{\Dx_k}\left(y_k,t\right)-u_{\Dx_\ell}\left(y_k,t\right)}+
\abs{u_{\Dx_\ell}\left(y_k,t\right)-u_{\Dx_\ell}\left(y_\ell,t\right)}.
$$
Both  terms on the right vanish as $k$ and $\ell$ become large since $\udx$ is
 uniformly continuous. Hence for any choice of $\seq{\theta_j}$, 
\eqref{eq:tildelim} holds.
In particular, this implies the pointwise convergence 
\begin{equation}
  \label{eq:otherlimits}
  \begin{gathered}
    \lim_{\Dx\to 0} \sum_j u_{j\pm 1/2}(t) \cha{I_j}(x) =
    u(x,t)\quad\text{and}\\
    \lim_{\Dx\to 0} \sum_j \tc_j(t)
    \cha{I_j}(x) = \frac{c'(u(x,t))}{4c(u(x,t))}=:\tc(u(x,t)),
\end{gathered}
\end{equation}
uniformly on compacts.
 
Next, we collect (in three lemmas) some well-known results related to weak
convergence. Throughout the paper we use overbars to denote weak limits.

 
\begin{lemma}[{\cite{Feireisl:NS_Book}}]\label{lem:weakconvlemma} 
Let $O$ be a bounded open subset of $\R^M$, with $M\ge 1$.
  
Let $\Set{v_n}_{n\ge 1}$ be a sequence of measurable functions on
$O$ for which
$$
\sup_{n\ge 1} \int_O \Phi(\abs{v_n(y)})\, dy<\infty,
$$
for some given continuous function $\Phi\colon[0,\infty)\to [0,\infty)$. 
Then along a subsequence as $n\to\infty$
$$
\text{$g(v_n)\weakto \overline{g(v)}$ in $L^1(O)$}
$$
for all continuous functions $g:\R\to \R$ satisfying
$$
\lim_{\abs{v}\to \infty} \frac{\abs{g(v)}}{\Phi(\abs{v})}=0.
$$
Let $g\colon \R\to (-\infty,\infty]$ be a lower semicontinuous convex
function and $\Set{v_n}_{n\ge 1}$ a sequence of measurable functions
on $O$, for which
$$
\textit{$v_n\weakto v$ in $L^1(O)$, $g(v_n)\in L^1(O)$ for each
$n$, $g(v_n)\weakto \overline{g(v)}$ in $L^1(O)$}.
$$
Then
$$
\text{$g(v)\le \overline{g(v)}$ a.e.~on $O$.}
$$
Moreover, $g(v)\in L^1(O)$ and
$$
\int_O g(v)\dy \le \liminf_{n\to\infty} \int_O g(v_n)\dy.
$$
If, in addition, $g$ is strictly convex on an open interval $(a,b)\subset \R$ and
$$
\text{$g(v)=\overline{g(v)}$ a.e.~on $O$},
$$
then, passing to a subsequence if necessary,
$$
\text{$v_n(y)\to v(y)$ for a.e.~$y\in \Set{y\in O\mid v(y)\in (a,b)}$.}
$$
\end{lemma}


Let $X$ be a Banach space and denote by $X^\star$ its dual. 
The space $X^\star$ equipped with the weak-$\star$ topology is denoted 
by $X^\star_{\mathrm{weak}}$, while $X$ equipped with the weak topology is 
denoted by $X_{\mathrm{weak}}$. According to the Banach--Alaoglu theorem, 
any bounded ball in $X^\star$ is $\sigma(X^\star,X)$-compact. 
If $X$ separable, then the weak-$\star$ topology 
is metrizable on bounded sets in $X^\star$, and thus one can consider the 
metric space $C\left([0,T];X^\star_{\mathrm{weak}}\right)$ 
of functions $v:[0,T]\to X^\star$ that are continuous with respect to 
the weak topology. We have $v_n\to v$ in $C\left([0,T];X^\star_{\mathrm{weak}}\right)$ 
if $\langle v_n(t),\phi \rangle_{X^\star,X}\to \langle v(t),\phi \rangle_{X^\star,X}$ 
uniformly with respect to $t$, for any $\phi\in X$. The following theorem is a 
straightforward consequence of the abstract Arzel\`a--Ascoli theorem:

\begin{lemma}[{\cite{Feireisl:NS_Book}}]
 \label{lem:timecompactness}
 Let $X$ be a separable Banach space, and suppose $v_n\colon [0,T]\to
 X^\star$, $n=1,2,\dots$, is a sequence of measurable functions such
 that
 $$
 \norm{v_n}_{L^\infty([0,T];X^\star)}\le C,
 $$
 for some constant $C$ independent of $n$. Suppose the sequence
 $$
 [0,T]\ni t\mapsto \langle v_n(t),\Phi \rangle_{X^\star,X}, \quad
 n=1,2,\dots,
 $$
 is equi-continuous for every $\Phi$ that belongs to a dense subset
 of $X$.  Then $v_n$ belongs to
 $C\left([0,T];X^\star_{\mathrm{weak}}\right)$ for every
 $n=1,2,\dots$, and there exists a $v\in
 C\left([0,T];X^\star_{\mathrm{weak}}\right)$ such that along a
 subsequence as $n\to \infty$
 $$
 \text{$v_n\to v$ in
   $C\left([0,T];X^\star_{\mathrm{weak}}\right)$}.
 $$
\end{lemma}

\begin{lemma}[Div-curl lemma {\cite{Maleketal}}]\label{lem:divcurl}
Let $Q\subset \R^2$ be a bounded domain.
Suppose
 \begin{align*}
      &v_\eps^1 \weakto \wlim{v}^1, \qquad v_\eps^2\weakto \wlim{v}^2,  \\
      &w_\eps^1 \weakto \wlim{w}^1, \qquad w_\eps^2\weakto \wlim{w}^2,
\end{align*}
in $L^2(Q)$ as $\eps\to 0$. Suppose also that the two 
sequences $\{\Div\left(v_\eps^1,v_\eps^2\right)\}_{\eps>0}$ and 
$\{ \Curl\left(w_\eps^1,w_\eps^2\right)\}_{\eps>0}$ 
lie in a common compact subset of $\Hneg(Q)$, where $\Div\left(v_\eps^1,v_\eps^2\right) 
=\partial_{x_1}v_\eps^1+\partial_{x_2}v_\eps^2$ 
and $\Curl\left(w_\eps^1,w_\eps^2\right) 
=\partial_{x_1}w_\eps^2-\partial_{x_2}w_\eps^1$.
Then along a subsequence
$$
\left(v_\eps^1,v_\eps^2\right)\cdot \left(w_\eps^1,w_\eps^2\right)
\to \left(\wlim{v}^1,\wlim{v}^2\right)\cdot
\left(\wlim{w}^1,\wlim{w}^2\right)\,\, 
\text{in $\mathcal{D}'(Q)$ as $\eps\to 0$.}
$$
\end{lemma}

\begin{lemma}[Murat's lemma \cite{Maleketal}]\label{Murat}
Suppose that $\left\{\Le\right\}_{\eps>0}$ is bounded in $\Sobolev{-1}{\infty}(\Pt)$. 
Suppose also that $\Le=\Le^1 + \Le^2$, where 
$\left\{\Le^1\right\}_{\eps>0}$ lies in a compact 
subset of $\Hneg(\Pt)$ and $\left\{\Le^2\right\}_{\eps>0}$ lies in a 
bounded subset of $\Mloc(\Pt)$. Then $\left\{\Le\right\}_{\eps>0}$ lies in 
a compact subset of $\Hneg(\Pt)$.
\end{lemma} 
According to Hypothesis \ref{hyp:main}, $R_0,S_0\in L^1(\R)\cap
L^p(\R)$ (with $p>3$). In view of \eqref{eq:l2bound},
\eqref{eq:lpbounds} and Lemma~\ref{lem:weakconvlemma}, there exist
$R,S\in L^\infty(0,T;L^q(\R))$, $q\in [1,p]$, $\wlim{R^2},\wlim{S^2}
\in L^\infty(0,T;L^r(\R))$, $r\in [1,p/2]$, such that along a
subsequence as $\Dx\to 0$
\begin{equation}
  \label{eq:weakconv-I}
  \begin{split}
    &\text{$R_{\Dx}\weakstar R$ in $L^{\infty}(0,T;L^{2}(\R))$},
    \quad 
    \text{$S_{\Dx}\weakstar S$ in $L^{\infty}(0,T;L^{2}(\R))$}, \\
    & \text{$R_{\Dx}^2 \weak \wlim{R^2}$ in $L^{\infty}(0,T;L^{r}(\R))$}, 
    \quad
    \text{$S_{\Dx}^2 \weakstar \wlim{S^2}$ in $L^{\infty}(0,T;L^{r}(\R))$}.\\
  \end{split}
\end{equation}
As a matter of fact, we can assume that for any 
function $f\in C^1(\R)$, with
\begin{equation}\label{eq:f-growth}
        \text{$\abs{f(z)}\le C\left(1+\abs{z}^2\right)$ and 
        $\abs{f'(z)}\le C\left(1+\abs{z}\right)$,}
\end{equation}
the following statements hold
\begin{equation}
        \label{eq:weak-f-limits}
        f(R_{\Dx}) \weakstar \wlim{f(R)},\quad  
        f(S_{\Dx}) \weakstar \wlim{f(S)}
        \quad \text{in $L^{\infty}(0,T;L^{p/2}(\R))$},
\end{equation} 
where the same subsequence of $\Dx\to 0$ applies to any $f$ from 
the specified class. Clearly, we can also assume that as $\Dx\to 0$
\begin{equation}\label{eq:svak2}
  \tc(\udx)\left(\Rdx-\Sdx\right)^2 \weakto \wlim{\tc(u)(R-S)^2} 
  =\tc(u)\wlim{(R-S)^2}\quad\text{in $L^{p/2}(\Pt)$,} 
\end{equation}
by equation \eqref{eq:otherlimits}.
From Lemmas~\ref{lem:Hlem} and \ref{lem:lpbounds}, it is not difficult 
to deduce that the functions
\begin{equation}
  \label{eq:weakcont-q}
  t \mapsto \int_{\R}  f(R_{\Dx}) \Phi \dx,\quad
  t \mapsto \int_{\R}  f(S_{\Dx}) \Phi \dx
\end{equation}
are equi-continuous on $[0,T]$ for every $\Phi\in C^\infty_0(\R)$.
In addition, $f(R_{\Dx})$ and $f(S_{\Dx})$ are bounded in
$L^{\infty}(0,T;L^{r}(\R))$, independently of $\Dx$. Now we apply
Lemma \ref{lem:timecompactness} with $X^\star=L^{r}(\R)$,
$X=L^{r'}(\R)$, and $r'=r/(r-1)$. Since $C^\infty_0(\R)$ is dense in
$L^r(\R)$, we can thus assume that $\wlim{f(R)},\wlim{f(S)}\in
C\left([0,T];L^{r}_{\mathrm{weak}}(\R)\right)$ and
\begin{equation}
        \label{eq:weak-f-limits-cont}
        \text{$f(R_{\Dx})\to \wlim{f(R)}$, $f(S_{\Dx})\to \wlim{f(S)}$ 
        in $C\left([0,T];L^{r}_{\mathrm{weak}}(\R)\right)$}.
\end{equation}
Of course, when $f(z)=z$, we can assume $\wlim{R},\wlim{S}\in C\left([0,T];L^{2}_{\mathrm{weak}}(\R)\right)$ and 
\begin{equation}
        \label{eq:weak-limits-cont}
        \text{$R_{\Dx}\to \wlim{R}$, $S_{\Dx}\to \wlim{S}$ 
        in $C\left([0,T];L^{2}_{\mathrm{weak}}(\R)\right)$.}
\end{equation}

\begin{lemma}\label{lem:Requation}
Assume Hypothesis \ref{hyp:main}. Then we have, cf.~\eqref{eq:wavesysCONS},
\begin{equation}
    \label{eq:Requation}
    \wlim{R}_t - \left(c(u)\wlim{R}\right)_x = -\tc(u) \,\wlim{(R-S)^2},
\end{equation}
and
\begin{equation}
    \label{eq:Sequation}
    \wlim{S}_t + \left(c(u)\wlim{S}\right)_x = -\tc(u)\,\wlim{(R-S)^2},
\end{equation}
in the sense of distributions on $\R\times [0,T)$, i.e., for any $\test \in C^\infty_0(\R\times [0,T))$,
\begin{multline*}
         \iint_{\Pt} \left(\wlim{R}\test_t - \left(c(u)\wlim{R}\right)\test_x \right) \,dx \,dt
        + \int_{\R} R_0(x) \test(x,0)\dx
        \\ = \iint_{\Pt} \tc(u) \,\wlim{(R-S)^2}\test\,dx \,dt
\end{multline*}
and
\begin{multline*}
        \iint_{\Pt} \left(\wlim{S}\test_t + \left(c(u)\wlim{S}\right)\test_x \right) \,dx \,dt
        + \int_{\R} S_0(x) \test(x,0)\dx
        \\  = \iint_{\Pt} \tc(u) \,\wlim{(R-S)^2}\,dx \,dt.
\end{multline*}
\end{lemma}

\begin{proof}
Fix $\test \in C^\infty_0(\R\times [0,T))$. The equation \eqref{eq:fRcons} 
with $f(R)=R$ reads
\begin{equation}\label{eq:xx}
    \frac{d}{dt}R_j - \Dp\left(\cj{-}R_j\right) = -\tc_j
    \left(R_j-S_j\right)^2.
\end{equation}
Set
\begin{equation*}
    \test_j(t)=\frac{1}{\Dx}\int_{I_j} \test(y,t)\,dy.
\end{equation*}
Next multiply the equation \eqref{eq:xx} with $\test_j$, sum over $j$, do a partial 
summation, integrate over $t$, to end up with
\begin{align}
  \label{eq:weakR}
  -\iint_{\Pt} \big(R_{\Dx} \test_t - c_{\Dx} R_{\Dx} &\test_x\big)
  \,dx\,dt + \int_\R R_{0,\Dx}(x)\test(x,0)\,dx \\
  &\quad=
  -\iint_{\Pt} \tc_{\Dx} (R_{\Dx} - S_{\Dx})^2 \test \,dxdt\notag\\
  &\qquad + \int_0^T \sum_j c_j R_j \int_{I_j} \left(\Dm \test_j -
    \test_x\right) \,dx \,dt\notag
\end{align}
where we have defined the functions  $c_{\Dx}$ and $\tc_{\Dx}$ by
\begin{equation*}
  c_{\Dx}(x,t)=\sum_j \cj{-}(t)\cha{I_{j-1/2}}(x)
  \quad\text{and}\quad \tc_{\Dx}(x,t)=\sum_j \tc_j(t)\cha{I_j}(x).
\end{equation*} 
By \eqref{eq:otherlimits}
\begin{equation}\label{eq:cdx-tcdx-conv}
  \text{$c_{\Dx}\to c(u)$, $\tc_{\Dx}\to \tc(u)$ 
    uniformly on $\mathrm{supp}\,(\test)$.}
\end{equation} 
Now
\begin{equation*}
  \abs{\Dm\test_j(t) - \test_x(x,t)} \le
  \norm{\test_{xx}}_{L^\infty(\R\times[0,T])} \Dx, \qquad 
  x\in I_j.
\end{equation*}
In view of this and \eqref{eq:lpbounds}, the last term in
\eqref{eq:weakR} is bounded as follows:
$$
\abs{ \int_0^T \sum_j c_j R_j \int_{I_j} \left(\Dm \test_j -
        \test_x\right) \,dx \,dt}\le C \Dx\to 0.
$$

Furthermore, in view of \eqref{eq:start-sterkkonv}, as $\Dx\to 0$
$$
\int_\R R_{0,\Dx}(x)\test(x,0)\,dx \to \int_\R R_{0}(x)\test(x,0)\,dx. 
$$

Hence, sending $\Dx \to 0$ in \eqref{eq:weakR} yields \eqref{eq:Requation}.
The evolution equation \eqref{eq:Sequation} 
for $\wlim{S}$ is proved in the same way. 
\end{proof}

We can also prove a generalization of the previous lemma.
\begin{lemma}\label{lem:frenormlemma}
Assume Hypothesis~\ref{hyp:main}, and let $f\in C^2(\R)$ be a convex 
function satisfying \eqref{eq:f-growth}. Then
\begin{align}
        \wlim{f(R)}_t - \left(c(u)\wlim{f(R)}\right)_x & \le
        2\tc(u)\left(\frac{1}{2}
        \wlim{f'(R)(R^2-S^2)}
        -\wlim{f(R)(R-S)}\right),\label{eq:flimren}\\
        \wlim{f(S)}_t + \left(c(u)\wlim{f(S)}\right)_x & \le
        -2\tc(u)\left(\frac{1}{2}\wlim{f'(S)(R^2-S^2)}
        - \wlim{f(S)(R-S)}\right),\label{eq:flimsren}
\end{align}
in the sense of distributions on $\R\times [0,T)$, i.e., for any
$\test \in C^\infty_0(\R\times [0,T))$, $\test\ge 0$,
\begin{align*}
  &\iint_{\Pt} \left(\wlim{f(R)}\test_t -
    \left(c(u)\wlim{f(R)}\right)\test_x \right) \,dx \,dt 
  + \int_{\R} f(R_0(x)) \test(x,0)\dx
  \\ & \quad  \ge - \iint_{\Pt}  2\tc(u)\left(\frac{1}{2}
    \wlim{f'(R)(R^2-S^2)}
    -\wlim{f(R)(R-S)}\right)\test\,dx \,dt
\end{align*}
and
\begin{align*}
  &\iint_{\Pt} \left(\wlim{f(S)}\test_t +
    \left(c(u)\wlim{f(S)}\right)\test_x \right) \,dx \,dt + \int_{\R}
  f(S_0(x)) \test(x,0)\dx \\ & \quad \ge \iint_{\Pt}
  2\tc(u)\left(\frac{1}{2}\wlim{f'(S)(R^2-S^2)}
    -\wlim{f(S)(R-S)}\right)\test\,dx \,dt.
\end{align*}
\end{lemma}
\begin{proof}
Similar to the proof of Lemma~\ref{lem:Requation}, starting 
from \eqref{eq:fRcons} and \eqref{eq:fScons}.
\end{proof}
The weak limits $\wlim{R^2}, \wlim{S^2}$ 
satisfy the initial data in a strong sense:
\begin{lemma}\label{lem:initsquare} 
Assume Hypothesis~\ref{hyp:main}. Then
\begin{equation}\label{eq:RS-squareinit}
  \begin{split}
    \lim_{t\to 0+} \int_{\R} \wlim{R^2}\,dx & 
    = \lim_{t\to 0+} \int_{\R} \wlim{R}^2\,dx
    =\int_\R R_0^2\,dx\\
    \lim_{t\to 0+} \int_{\R} \wlim{S^2}\,dx &
    =\lim_{t\to 0+} \int_{\R} \wlim{S}^2\,dx
    =\int_\R S_0^2\,dx.
  \end{split}
\end{equation}
\end{lemma}
\begin{proof} 
  Since $\wlim{R},\wlim{S}\in
  C\left([0,T];L^{2}_{\mathrm{weak}}(\R)\right)$, it follows from
  \eqref{eq:Requation}, \eqref{eq:Sequation} that
  $$
  \text{$\wlim{R}(\cdot,t) \weak R_0$, $\wlim{S}(\cdot,t) \weak S_0$ 
    in $L^{2}(\R)$ as $ t\to 0+$.}
  $$
  From this, \eqref{eq:weakconv-I}, and Lemma
  \ref{lem:weakconvlemma} we conclude that
  \begin{equation}
    \label{eq:timecont-tmp-I}
    \int_\R R_0^2\dx\le \liminf_{t\to 0+}\int_\R \wlim{R}^2\dx,
    \quad
    \int_\R S_0^2\dx \le\liminf_{t\to 0+}\int_\R \wlim{S}^2\dx.
  \end{equation}
  On the other hand, \eqref{eq:weakconv-I} says that $R_{\Dx}(\cdot,t)
  \weak \wlim{R}(\cdot,t)$, $S_{\Dx}(\cdot,t) \weak \wlim{S}(\cdot,t)$
  in $L^2(\R)$ for a.e.~$t>0$, and thereby, using also
  \eqref{eq:l2bound} and \eqref{eq:start-sterkkonv},
  \begin{equation}
    \label{eq:time-cont-tmp-II}
    \int_\R \left(\wlim{R}^2+\wlim{S}^2\right)(t,x) \dx\le 
    \int_\R \left( \wlim{R^2}+\wlim{S^2}\right)(t,x) \dx 
    \le \int_\R \left(R_0^2+S_0^2\right)\dx.
  \end{equation}
  Since $\wlim{R^2}, \wlim{S^2}\in
  C\left([0,T];L^{r}_{\mathrm{weak}}(\R)\right)$ (with $r>1$), one can
  prove that this inequality actually holds for all $t>0$. Combining
  \eqref{eq:timecont-tmp-I} and \eqref{eq:time-cont-tmp-II} yields
  \eqref{eq:RS-squareinit}.
\end{proof}


\begin{lemma}\label{lem:weakconvRS} 
Assume Hypothesis \ref{hyp:main}, and let and let $f,g\in C^2(\R)$ be 
functions satisfying $\abs{f(z)},\abs{g(z)}\le C \abs{z}$. Then 
\begin{equation}\label{eq:prod}
        f(\Rdx)g(\Sdx) \to \wlim{f(R)}\, \wlim{g(S)}
        \quad\text{in the distributional sense on $\R\times (0,T)$.}
\end{equation}
\end{lemma}

\begin{proof}
We will show that the sequences
\begin{align*}
        & \seq{\partial_t f(\Rdx) 
        -\partial_x\left(c\left(\udx\right)f(\Rdx)\right)}_{\Dx>0}, \\
        & \seq{\partial_t g(\Sdx) 
        +\partial_x\left(c\left(\udx\right)g(\Sdx)\right)}_{\Dx>0}
\end{align*}
are compact in $\Hneg(\R\times (0,T))$. 

Introducing the distribution $\mathcal{L}_{\Dx}=\partial_t
f(\Rdx)-\partial_x\left(c\left(\udx\right)f(\Rdx)\right)$, we find
\begin{align}
  & \langle \mathcal{L}_{\Dx},\test\rangle \notag
  \\ & =-\iint_{\Pt} \left[2\tc_{\Dx}
    \left(\frac{1}{2}f'(\Rdx)\left(\Rdx^2-\Sdx^2\right)-f(\Rdx)(\Rdx-\Sdx\right)
    +C_f\right]\test  \,dxdt\label{eq:firstterm} \\
  &\qquad + \int_0^T \sum_j c_j f\left(R_j\right)
  \int_{I_j} \left(\Dm 
    \test_j - \test_x\right) \,dx\,dt,\label{eq:secondterm}
\end{align}
for $\phi\in C^\infty_0(\R\times (0,T))$, where $C_f(x,t)$ is a function 
that is bounded in $L^1(\R\times (0,T))$ independently of 
$\Dx$, cf.~ \eqref{eq:weakR}. The last term above is bounded by 
$$
p_{\Dx}:=C\norm{R_0}_{L^2(\R)}\int_0^T\bigl\|\sum_j
  \Dm\test_j\cha{I_j}-\test_x\bigr\|_{L^2(\R)}\,dt.
$$
Since $\sum_j\Dm\test_j\cha{I_j}$ is a piecewise constant approximation to
$\test_x$, by Lemma~\ref{lem:constl2converg}, $p_{\Dx}$
tends to zero as $\Dx\to 0$. Thus
 we infer that
\begin{equation*}
  \abs{\langle \mathcal{L}_{\Dx},\test\rangle}
  \le  C\norm{\test}_{L^\infty(\R\times (0,T))}+p_{\Dx},
\end{equation*}
where $p_{\Dx}$ tends to zero with $\Dx$.  Thus \eqref{eq:secondterm} is in
a compact subset of $\Hneg(\R\times (0,T))$, while \eqref{eq:firstterm} is
in a bounded subset of the locally bounded Radon measures.  Hence
Murat's lemma implies that $\mathcal{L}_{\Dx}$ is compact in
$\Hneg(\R\times (0,T))$.  By analogous arguments, $\seq{\partial_t
  g(\Sdx)+\partial_x\left(c\left(\udx\right)g(\Sdx)\right)}$ is
compact in $\Hneg(\R\times (0,T))$.

Now by the div-curl lemma on the sequences  
$$
\seq{g(\Sdx),c\left(\udx\right)g(\Sdx)}_{\Dx>0}\quad\text{and}\quad
\seq{c\left(\udx\right)f(\Rdx),f(\Rdx)}_{\Dx>0}, 
$$
we see that
\begin{align}
  \label{eq:dclem}
  2c\left(\udx\right)f(\Rdx)g(\Sdx) 
  &\to 2c(u) \wlim{f(R)}\,\wlim{g(S)} 
  \quad \text{in the distributional sense,}
\end{align}
which, due to \eqref{eq:dclem} and Lemma \ref{lem:udxconverg}, concludes 
the proof of \eqref{eq:prod}.
\end{proof}

An immediate consequence of the previous lemma is the following result.
\begin{corollary}
There holds
\begin{equation}
        \label{eq:squarelim}
        \wlim{(R-S)^2} = \wlim{R^2} - 2\wlim{R}\,\wlim{S} + \wlim{S^2} 
        \quad \text{a.e.~in $\R\times (0,T)$.}
\end{equation}
\end{corollary}

\begin{proof}
Since we can assume that $\Rdx \Sdx\weakto \wlim{RS}$ in $L^1(\R\times (0,T))$, it follows 
from Lemma \ref{lem:weakconvRS} that 
$$
\iint_{\R\times (0,T)}  \wlim{RS}\test \,dxdt=
\iint_{\R\times (0,T)}\wlim{R}\, \wlim{S} \test \,dxdt, 
\quad \forall \test\in C^\infty_0(\R\times (0,T)),
$$
from which we infer that $\wlim{RS}=\wlim{R}\, \wlim{S}$ a.e.; Hence 
\eqref{eq:squarelim} follows. 
\end{proof}

\begin{lemma}\label{lem:strongconv} 
Assume Hypothesis~\ref{hyp:main}. Then 
\begin{equation}
        \label{eq:strongconv}
        \wlim{R^2}=\wlim{R}^2\quad\text{and}\quad \wlim{S^2}=\wlim{S}^2,
        \quad \text{for a.e.~$(x,t)\in\Pt$.}
\end{equation}
Consequently, as $\Dx\to 0$,
\begin{equation}
        \label{eq:strongconv1}
        \text{$\Rdx \to \wlim{R}$, $\Sdx\to \wlim{S}$ 
        in $L^{2}_{\loc}(\Pt)$ and almost everywhere in $\Pt$}.
\end{equation}

\end{lemma}

\begin{proof}
Using Lemmas \ref{lem:Requation}--\ref{lem:initsquare} and Corollary \ref{eq:squarelim}, we can argue 
exactly as in, e.g., Zhang and Zheng \cite{ZZ:2004}, to arrive at \eqref{eq:strongconv} and 
\eqref{eq:strongconv1}.
\end{proof}

\begin{lemma}\label{lem:ueqn1}
Assume Hypothesis \ref{hyp:main}. 
 Then $u$ is a weak solution of \eqref{eq:waveeq}, i.e.,
  \begin{equation*}
    \frac{\partial^2 u}{\partial t^2}-c(u)\frac{\partial}{\partial
      x}\left( c(u) \frac{\partial u}{\partial x}\right) = 0
  \end{equation*}
  weakly in $\Pt$ in the sense that
\begin{equation}\label{eq:usubsol}
  \iint_{\Pt}\big(u_t\test_t -c(u)_x(c(u)\test)_x\big)\, dxdt 
  = 0
\end{equation} 
for all test functions $\test\in C^\infty_0(\Pt)$. Here $u_t$ and $c(u)_x$ 
are given by \eqref{eq:ut} and \eqref{eq:dclim},
respectively.
\end{lemma}
\begin{proof}

We claim that
\begin{equation}
        \label{eq:dclim}
        c(u)_x = 2\tc(u)(\wlim{R}-\wlim{S}),\quad\text{weakly.}
\end{equation} 

To this end let
\begin{equation}\label{eq:chdef}
        c_{\Dx}=\sum_j \cj{-}\cha{I_j},
\end{equation}
and compute 
\begin{align*}
        \left\langle \frac{\partial}{\partial x}c_{\Dx},\test\right\rangle
        &= -\iint_{\Pt} c_{\Dx}\test_x\,dxdt\\
        &=-\int_0^T \sum_j \int_{I_j} \cj{-} \test_x \, dxdt\\
        &=-\int_0^T \sum_j \cj{-} \Dp\test\left(x_{j-1/2},t\right)\Dx\,dt\\
        &=\int_0^T \sum_j \left(\Dp\cj{-}\right)
        \test\left(x_{j-1/2},t\right)\Dx \,dt\\
        &=\int_{0}^T \Dx\sum_j \test\left(x_{j-1/2},t\right) 
        \left[2\tc_j \left(R_j-S_j\right)\right]\,dt.
\end{align*}
By sending $\Dx$ to zero in this equality, and using 
\eqref{eq:otherlimits}, our claim \eqref{eq:dclim} follows.

From Lemma \ref{lem:Requation}, we find that
\begin{equation} \label{eq:weakRS}
        \left(\wlim{R}-\wlim{S}\right)_t 
        -\left(c(u)(\wlim{R}+\wlim{S})\right)_x = 0, 
        \quad \text{in the sense of distributions.}
\end{equation}
Observe that for a function $u$ that is at least one time differentiable we have that 
$(c(u)u_t)_x=(c(u)u_x)_t$ holds in the distributional sense. Specifically, we have
\begin{align}
        \iint_{\Pt} (c(u) u_t)\test_x \,dx dt&= \iint_{\Pt} c(u)\test_x u_t\,dx dt
        =\iint_{\Pt} \big((c(u)\test)_x- c'(u)u_x\test \big) u_t \,dx dt\notag \\
        &= \iint  \big((c(u)\test)_t - c'(u)u_t\test \big) u_x\,dx dt \notag \\
        & = \iint (c(u) u_x)\test_t \,dx dt.
\end{align} 
Thus we see that this can be rewritten
\begin{equation*}
\frac{\partial}{\partial x}\left( c(u)\left(2u_t -
    (\wlim{R}+\wlim{S})\right)\right) =0
\end{equation*}
in the sense of distributions. Hence
\begin{equation} \label{eq:ut}
        u_t = \frac{1}{2}(\wlim{R}+\wlim{S}),
\end{equation}
since 
\begin{equation*}
        \lim_{x\to -\infty}
        u(x,t)=\lim_{x\to-\infty}(\wlim{R}(x,t)
        +\wlim{S}(x,t)) = 0.
\end{equation*}

Set
\begin{equation*}
        \wlim{R}^\eps(x,t)=\int_\R \wlim{R}(y,t)j^\eps(x-y)\,dy,
\end{equation*}
where $j^\eps$ is a standard mollifier. Then
\begin{equation*}
        \wlim{R}^\eps_t - (c(u) \wlim{R}^{\eps})_x = -\tc(u)\wlim{(R-S)^2}*j^\eps + r^\eps,
\end{equation*}
where by the DiPerna--Lions folklore lemma 
\begin{equation*}
        r^\eps=(c(u)\wlim{R})_x*j^{\eps} - (c(u) \wlim{R}^{\eps})_x 
\end{equation*}
tends to zero in $L^1_{\mathrm{loc}}(\Pt)$. This in turn implies that
\begin{align}
        \wlim{R}^\eps_t - c(u) \wlim{R}^\eps_x &= -\tc(u)\wlim{(R-S)^2}*j^\eps 
        +2\tc(u)(\wlim{R}-\wlim{S}) \wlim{R}^\eps + r^\eps\notag\\
        &=-\tc(u) \left(\wlim{(R-S)^2} 
        - 2 \wlim{R}^{\eps}(\wlim{R}-\wlim{S})\right) + r^\eps.\label{eq:Repsnoncons}
\end{align}
Similarly, with $ \wlim{S}^{\eps}=\wlim{S}*j^{\eps}$, we get
\begin{equation}\label{eq:Sepsnoncons}
        \wlim{S}^\eps_t + c(u) \wlim{S}^\eps_x = 
        -\tc(u) \left(\wlim{(R-S)^2} + 2\wlim{S}^\eps (\wlim{R}-\wlim{S})\right) 
        - s^{\eps},
\end{equation}
where $s^{\eps}$ tends to zero in $L^1_{\mathrm{loc}}(\Pt)$. 
Adding \eqref{eq:Repsnoncons} and \eqref{eq:Sepsnoncons} and sending 
$\eps$ to zero, we get, after using Lemma~\ref{lem:strongconv},
\begin{equation*}
        \frac{\partial}{\partial t}\frac{1}{2}(\wlim{R}+\wlim{S}) 
        -c(u)\frac{\partial}{\partial x}\frac{1}{2} (\wlim{R}-\wlim{S})
        =\tc(u)\left(\left(\wlim{R}-\wlim{S}\right)^2 - \wlim{(R-S)^2}\right)=0,
\end{equation*}
which, using \eqref{eq:ut}, \eqref{eq:dclim} and Lemma~\ref{lem:weakconvRS}, can 
be rewritten as \eqref{eq:usubsol}.
\end{proof}
We collect some of our results in the following theorem.
\begin{theorem}
  \label{thm:l3negconv} 
  Assume \eqref{c_est} and Hypothesis \ref{hyp:main}. Then the
  semi-discrete difference scheme defined by
  \eqref{eq:Reqn}--\eqref{eq:tcdef} produces a sequence that converges
  to a weak solution of \eqref{eq:waveeq}.

The same conclusion holds if Hypothesis \ref{hyp:main} is replaced by the assumption that 
$u_0\in W^{1,1}(\R)$ and $v_0\in L^{1}(\R)$ and that $R_0$ and $S_0$ take values in $[-M,0]$ for some positive
  constant $M$.
 \end{theorem}
 \begin{proof}
 Observe that in the case when $R_0$ and $S_0$ take values in $[-M,0]$, Lemma \ref{lem:invariant} shows that  $R_j(t), S_j(t)\in[-M,0]$ for all $t$. Thus $\Rdx,\Sdx\in L^\infty\cap L^1$, and by interpolation we see that Hypothesis \ref{hyp:main} is satisfied.
 \end{proof}
\begin{remark}\label{rem:anotherscheme}
  \normalfont The equation \eqref{eq:ujdef2} is awkward to use in
  practice, since one must compute the inverse of $F$ at
  each node. In order to circumvent this we may redefine $\tc_j$ slightly.
  Let $R_j$ and $S_j$ be defined as before, but let $u_{j-1/2}$ be
  defined by
  \begin{equation}
    \label{eq:ujm1def}
    \Dp \uj{-} = \frac{R_j - S_j}{c(\uj{-})+c(\uj{+})}.
  \end{equation}
  This is also a nonlinear equation to solve for $\uj{+}$, but solving
  this is likely to be easier than inverting $F$. However, if we
  accept a certain imbalance, we can define $\uj{+}$ by
  \begin{equation}
    \label{eq:ujm2def}
    \Dp \uj{-} = \frac{R_j - S_j}{2c(\uj{-})}.
  \end{equation}
  In order to get our approach to work, it is essential that
  \eqref{eq:dccrucial} holds. Therefore, we shall define $\tc_j$ so
  that this is the case. Since $c$ is continuous, we have that
  \begin{equation*}
  \Dp c(\uj{-}) = c'\left(\bar{u}_j\right) \Dp \uj{-},
  \end{equation*}
  for some $\bar{u}_j$ between $\uj{-}$ and $\uj{+}$.  If $\uj{+}$
  is defined by \eqref{eq:ujm1def} then we set
  \begin{equation}
    \label{eq:tc1def}
    \tc_j = \frac{c'(\bar{u}_j)}{2\left(c(\uj{-})+c(\uj{+})\right)},
  \end{equation}
  while if $\uj{+}$ is defined by \eqref{eq:ujm2def} we set
  \begin{equation}\label{eq:tc2def}
    \tc_j = \frac{c'(\bar{u}_j)}{4c\left(\uj{-}\right)}.
  \end{equation}
  In both cases \eqref{eq:dccrucial} holds. Therefore, the schemes
  defined by \eqref{eq:Reqn}, \eqref{eq:Reqn}, and \eqref{eq:ujm1def},
  \eqref{eq:tc1def} or \eqref{eq:ujm2def}, \eqref{eq:tc2def} all
  produce sequences converging to weak solutions of \eqref{eq:waveeq}.
\end{remark}

\section{Numerical examples}\label{sec:numer}
The semi-discrete scheme defined by \eqref{eq:Reqn}--\eqref{eq:initapprox} is rather involved, in particular the computation of the quantity $\tilde c_j$. For actual computations one needs to make a further discretization of the time variation. 
We have considered the following versions of the semi-discrete scheme
defined by \eqref{eq:Reqn}--\eqref{eq:dccrucial}. These schemes all
use an explicit discretization of \eqref{eq:Reqn}--\eqref{eq:Seqn},
\begin{equation}
  \label{eq:RSexplicit}
  \begin{aligned}
    \Dtp R^n_j - c^n_{j+1/2}\Dp R^n_j &= \tc^n_j\left((R^n_j)^2 -
      (S^n_j)^2\right), \\
    \Dtp S^n_j + c^n_{j-1/2}\Dm S^n_j &= -\tc^n_j\left((R^n_j)^2 -
      (S^n_j)^2\right),
  \end{aligned}
\end{equation}
where
\begin{equation*}
c^{n}_{j-1/2} = c\left(u^n_{j-1/2}\right)
\end{equation*}
and
\begin{equation*}
\Dtp K^n_j = \frac{1}{\Dt}\left(K^{n+1}_j - K^{n}_j\right).
\end{equation*}
Furthermore, since we wish something like \eqref{eq:dccrucial} to
hold,
\begin{equation}\label{eq:dc_discrete}
  \tc^n_j = \frac{\Dp c^n_{j-1/2}}{2 (R^n_j -S^n_j)}.
\end{equation}
The difference between these schemes consists in the way the
``cofficients'' $c^n_{j-1/2}$ are defined.
\begin{enumerate}
\item {\bf Integration in time.}  We update $u^n_{j-1/2}$ by
  considering a discrete version of \eqref{eq:ut},
  \begin{multline}
    \label{eq:u_timedef}
    u^{n+1}_{j+1/2} = u^{n}_{j+1/2} \\
    + \frac{\Dt}{8}\left(
      R^{n}_{j}+R^{n}_{j+1} +  R^{n+1}_{j}+R^{n+1}_{j+1} + 
      S^{n}_{j}+S^{n}_{j+1} +  S^{n+1}_{j}+S^{n+1}_{j+1}\right).
  \end{multline}
  We use this scheme since $u$ is discretized on a grid that is
  staggered with respect to that of $R$ and $S$.
\item{\bf Integration in space.} Knowing $u^n_{-N-1/2}$ for some large
  $N$, we can set
  \begin{equation}
    \label{eq:duinspace}
    u^{n}_{j+1/2} = u^n_{j-1/2} + \Dx\frac{R^n_j -
    S^n_j}{c(u^n_{j-1/2})+c(u^n_{j+1/2})}.
  \end{equation}
\end{enumerate}
In this section we describe two examples. 

\begin{figure}
  \centering
  \begin{tabular}{c}
    \includegraphics[width=0.95\linewidth]{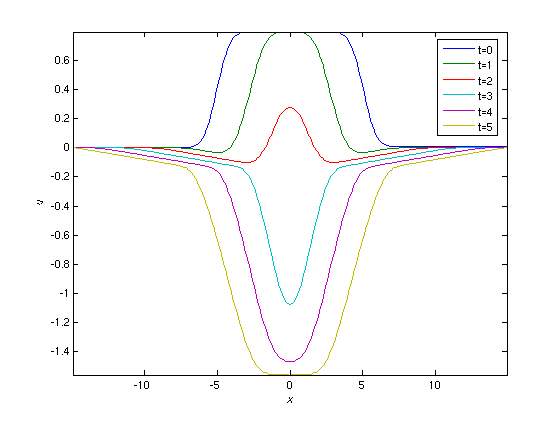} \\
    \includegraphics[width=0.95\linewidth]{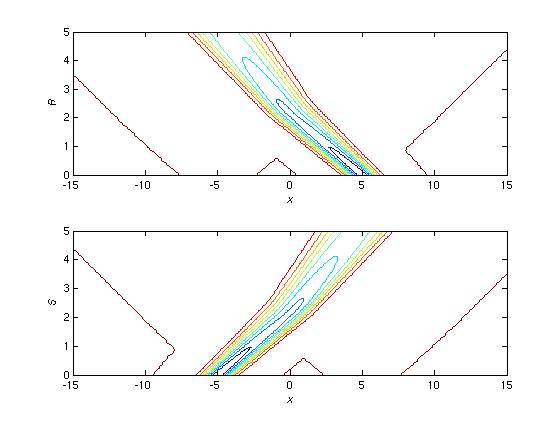}%
  \end{tabular}
  \caption{The scheme \eqref{eq:duinspace} with the initial data
    \eqref{eq:RSnegdef} and $c$ given by \eqref{eq:cincreasedef}. The
    $u$ variable (top), the $R$ variable (middle) and the $S$
    variable (bottom) as functions of $x$ and $t$.}
  \label{fig:2}
\end{figure}

\begin{figure}
   \centering
   \begin{tabular}{c}
     \includegraphics[width=0.9\linewidth]{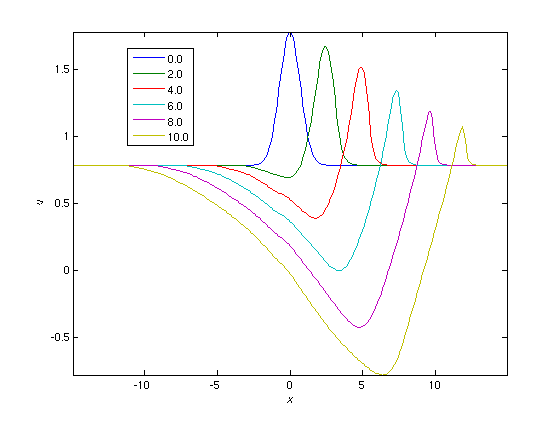} \\
     \includegraphics[width=0.9\linewidth]{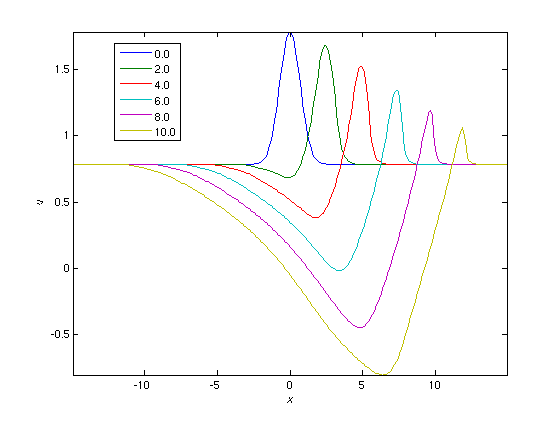}%
   \end{tabular}
   \caption{The scheme \eqref{eq:u_timedef} (top) and
     \eqref{eq:duinspace} (bottom), with the initial data
     \eqref{eq:initialdata} and $c$ given by \eqref{eq:cnumeric}.  } 
   \label{fig:1}
 \end{figure}

Consider first the case where the function $c$ is given by
\begin{equation}
  \label{eq:cincreasedef}
  c(u)=\frac{2}{\pi}\left(\pi + \arctan(u)\right),
\end{equation}
and the initial data are given by
\begin{equation}
  \label{eq:RSnegdef}   
  R(x,0)=-2e^{-(x-5)^2},\quad S(x,0)=2e^{-(x+5)^2}.
\end{equation}
In this case
\begin{equation*}
u(x,0)=\int_{-\infty}^x \frac{R(y,0)-S(y,0)}{2c(u(y,0))}\,dy,\quad
u_t(x,0)=\frac{1}{2}(R(x,0)+S(x,0)).
\end{equation*}
In Figure~\ref{fig:2} we show the computed solution $u$ (top) and
$R$ and $S$ (bottom).
The discrete difference scheme can be studied numerically also in
cases not covered by the convergence results in this paper.  We have
included an example of that type here.  Define the function $c$ by
\begin{equation}
  \label{eq:cnumeric}
  c(u)=\sqrt{\alpha \cos^2(u) + \beta \sin^2(u)},\quad
  \alpha=1.5,\quad \beta=0.5.  
\end{equation}
When testing, we take the initial data from
\cite{GlasseyHunterZheng:sing}, and use
\begin{equation}
  \label{eq:initialdata}
  u(x,0)=\frac{\pi}{4}+e^{-x^2}, \quad
  u_t(x,0)=-c(u(x,0))\frac{\partial }{\partial x} u(x,0).
\end{equation}
In order for the two schemes to be compatible, we have defined
$u^0_{j-1/2}$ by
\begin{equation*}
\Dp u^0_{j-1/2} =
\frac{R_j^0-S^0_j}{c(u^0_{j+1/2})+c(u^0_{j-1/2})},
\end{equation*}
even for the scheme using by \eqref{eq:u_timedef}. In
Figure~\ref{fig:1} we show $u$ for the two methods with initial data
\eqref{eq:initialdata} using $\Dx=30/256$, and $\Dt=\Dx$. We remark
that using $\Dt=\Dx/M$ where $M$ is a large integer, produced very
similar results.

\appendix
\section{Higher integrability properties}\label{sec:higher}
\renewcommand{\theequation}{A.\arabic{equation}}
\renewcommand{\thetheorem}{A.\arabic{theorem}} \setcounter{theorem}{0}
\setcounter{equation}{0}

In this appendix we prove a so-called \emph{higher integrability}
result. Briefly stated, we have that if $R_0$ and $S_0$ are nonpositive 
and in $L^1\cap L^2$, then $\partial_x u(\dott,t)$ is in
$L^{p}_{\mathrm{loc}}$ for all $p\in [2,3)$. This is obvious if $R_0$
and $S_0$ are in $L^3$, and the significance of this section is that
the $\partial_x u$ is more integrable than is to be expected. The
reason for including this is that we suspect that such a property will
play a role in a (yet unknown) uniqueness result.

Throughout the appendix we assume that $R_0$ and $S_0$ are nonpositive
and in $L^1(\R)\cap L^2(\R)$. Then by Lemma~\ref{lem:lpbounds} 
$$
\Dx\sum_j\left( \abs{R_j}^{1+\alpha}+\abs{S_j}^{1+\alpha} \right)\le C,
$$
for any $\alpha\in [0,1]$ and some constant $C$ which is independent
of $\Dx$. We also recall that for any smooth function $f$ we have that
\begin{equation}\label{eq:afnon}
\frac{d}{dt} f_j - \cjph\Dp f_j - \frac{\Dx}{2}f''_j
\left(\Dp R_j\right)^2
= \tc_j f'_j \left(R_j^2-S_j^2\right),
\end{equation}
where $f_j=f(R_j)$, $f'_j=f'(R_j)$ and $f''_j=f''(r_j)$ for some $r_j$
between $R_j$ and $R_{j+1}$. 

We now let $\alpha$ be a constant in $[0,1)$ and define $f$ to be a
$C^\infty$ function such that
\begin{align*}
f'(K)&=
\begin{cases}
  0, \quad &K>-1/2,\\
  \abs{K}^\alpha, & K<-1,
\end{cases}\\
f(K)&=\int_0^K f'(\sigma)\,d\sigma = 
\begin{cases}
  0, \quad & K>-1/2,\\
  \frac{-\abs{K}^{1+\alpha}}{1+\alpha} + C, &K<-1.
\end{cases}
\end{align*}
Note that $f''(K)$ is bounded.
Let $\chi(x)$ be a smooth function such that $0\le \chi(x)\le 1$ and
$$
\chi(x)=
\begin{cases}
  0, \quad &x\not\in [a-1,b+1],\\
  1, \quad &x\in [a,b],
\end{cases}
$$
where $a<b$ are real numbers. Set $\chi_j=\chi(x_j)$. We multiply
\eqref{eq:afnon} by $\chi_j\Dx$, sum over $j$ and integrate over
$[0,T]$ to end up with
\begin{align*}
  \Dx\sum_j f_j \chi_j \bigm|^T_0 -  \int_0^T \Dx \sum_j\chi_j \cjph
  \Dp f_j \,dt\; + &\int_0^T \Dx \sum_j\chi_j \frac{\Dx}{2} f''_j \left(\Dp
    R_j\right)^2 \,dt \\
  &= \int_0^T \Dx \sum_j \chi_j f'_j\tc_j\left(R_j^2 - S_j^2\right)\,dt. 
\end{align*}
After a partial summation of the second term on the right, we obtain
\begin{align*}
   \Dx\sum_j f_j \chi_j \bigm|^T_0  &+
   \int_0^T \Dx \sum_j f_j
   \left[\chi_j2\tc_j\left(R_j-S_j\right)+\cjmh \Dm \chi_j\right]
    \,dt\\ &- \int_0^T \Dx \sum_j \frac{\Dx}{2}\chi_j f''_j \left(\Dp
      R_j\right)^2 \,dt 
    = \int_0^T \Dx \sum_j \chi_j f'_j\tc_j\left(R_j^2 - S_j^2\right)\,dt. 
\end{align*}
Rearranging this we find that
\begin{align*}
  \int_0^T\Dx \sum_j \chi_j 2\tc_j &\left[
    \left(R_j-S_j\right)\left(\frac{-\abs{R_j}^{1+\alpha}}{1+\alpha}\right) - 
    \frac{1}{2}\left(R_j^2-S_j^2\right)R_j^\alpha\right]\cha{\{R_j<-1\}}\,dt\\
  &=- \Dx\sum_j f_j \chi_j \bigm|^T_0 -\; C \int_0^T\Dx \sum_j \chi_j 2\tc_j 
  \left(R_j-S_j\right)\,dt \\
  &\quad \int_0^T \Dx\sum_j \chi_j\tc_j  \left[
    \left(R_j^2-S_j^2\right)f'_j - 2\left(R_j-S_j\right)f_j\right]\cha{\{R_j>-1\}}
  \, dt\\ 
  &\quad - \int_0^T \Dx\sum_j \chi_j \frac{\Dx}{2} f''_j \left(\Dp
    R_j\right)^2 \,dt. 
\end{align*}
By the $L^p$ estimates, Lemma~\ref{lem:lpbounds}, all terms on the
right-hand side of this are bounded by a constant $C_{T,a,b}$
depending only on $T$, $a$, $b$ and on the $L^1$ and $L^2$ norms of
$R_0$ and $S_0$. Furthermore, by the same lemma, 
$$
\biggl| \int_0^T \Dx\sum_j   
\chi_j 2\tc_j \left[
    \left(R_j-S_j\right)\left(\frac{-\abs{R_j}^{1+\alpha}}{1+\alpha}\right) - 
    \frac{1}{2}\left(R_j^2-S_j^2\right)R_j^\alpha\right]\cha{\{R_j>-1\}}\,dt\biggr|
  \le C_{T,a,b}. 
$$ 
Therefore we get the bound
\begin{align}
  \biggl| \int_0^T\Dx\sum_j \chi_j \tc_j \left[
    \left(R_j-S_j\right)\left(\frac{-\abs{R_j}^{1+\alpha}}{1+\alpha}\right) - 
    \frac{1}{2}\left(R_j^2-S_j^2\right)R_j^\alpha\right]\,dt\biggr| &\le
  C_{T,a,b},\\
  \intertext{and similarly}
  \biggl|
  \int_0^T\Dx\sum_j \chi_j \tc_j \left[
    -\left(R_j-S_j\right)\left(\frac{-\abs{S_j}^{1+\alpha}}{1+\alpha}\right)
    -\frac{1}{2}\left(S_j^2-R_j^2\right)S_j^\alpha\right]\,dt
  \biggr|&\le C_{T,a,b}.
\end{align}
Adding these two and recalling that $\abs{R_j}=-R_j$ and
$\abs{S_j}=-S_j$, we get the bound
\begin{multline}\label{eq:m1}
  \biggl| \int_0^T\Dx\sum_j \chi_j \tc_j\biggl[\frac{1}{1+\alpha}\left(
    \abs{R_j}^{1+\alpha}-\abs{S_j}^{1+a}\right)\left(\abs{R_j}-\abs{S_j}\right)\\
  -
  \frac{1}{2}\left(\abs{R_j}^2-\abs{S_j}^2\right)
  \left(\abs{R_j}^\alpha-\abs{S_j}^\alpha\right)  \biggr] \,dt\biggr|\le
  C_{T,a,b}. 
\end{multline} 
The term in the square bracket above can be rewritten as 
\begin{align*}
  0\le \frac{1}{1+\alpha}&\left(
    \abs{R_j}^{1+\alpha}-\abs{S_j}^{1+a}\right)\left(\abs{R_j}-\abs{S_j}\right)\\
  &\qquad -
  \frac{1}{2}\left(\abs{R_j}^2-\abs{S_j}^2\right)
  \left(\abs{R_j}^\alpha-\abs{S_j}^\alpha\right) \\
  &= \left(\frac{1}{1+\alpha}-\frac{1}{2}\right) 
  \left(\abs{R_j}-\abs{S_j}\right)
  \left(\abs{R_j}^{1+\alpha}-\abs{S_j}^{1+\alpha}\right) \\
  &\qquad + \frac{1}{2}
  \abs{R_j}^\alpha\abs{S_j}^\alpha\left(\abs{R_j}-\abs{S_j}\right)
  \left(\abs{R_j}^{1-\alpha}-\abs{S_j}^{1-\alpha}\right)\\
  &=\frac{1}{2(1+\alpha)} 
  \biggl[
  (1-\alpha)\left(\abs{R_j}-\abs{S_j}\right) \\
 &\qquad \times \left(\abs{R_j}^{1+\alpha}-\abs{S_j}^{1+\alpha}
    +\abs{R_j}^\alpha\abs{S_j}^\alpha 
    \left(\abs{R_j}^{1-\alpha}-\abs{S_j}^{1-\alpha}\right)
  \right)
  \\
  &\hphantom{=\frac{1}{2(1+\alpha)} \biggl[}\quad
  +2\alpha \left(\abs{R_j}-\abs{S_j}\right)\abs{R_j}^\alpha\abs{S_j}^\alpha
  \left(\abs{R_j}^{1-\alpha}-\abs{S_j}^{1-\alpha}\right)\biggr]
  \\
  &=\frac{1-\alpha}{2(1+\alpha)}\left(\abs{R_j}-\abs{S_j}\right)^2
  \left(\abs{R_j}^\alpha+\abs{S_j}^\alpha\right)
  \\
  &\qquad + 
  \frac{\alpha}{1+\alpha} 
  \left(\abs{R_j}-\abs{S_j}\right)\abs{R_j}^\alpha\abs{S_j}^\alpha
  \left(\abs{R_j}^{1-\alpha}-\abs{S_j}^{1-\alpha}\right).  
\end{align*}
Hence, multiplying \eqref{eq:m1} by $2(1+\alpha)$, we arrive at
\begin{equation}
  \label{eq:intest_a1}
  \begin{aligned}
    \int_0^T\Dx\sum_j \chi_j \tc_j &\biggl[ (1-\alpha)\left(\abs{R_j}-\abs{S_j}\right)^2
    \left(\abs{R_j}^{\alpha}+\abs{S_j}^{\alpha}\right)\\
    &+2\alpha\abs{R_j}^\alpha\abs{S_j}^\alpha
      \left(\abs{R_j}^{1-\alpha}-\abs{S_j}^{1-\alpha}\right)\biggr]\,dt\le
    C_{T,a,b}.
  \end{aligned}
\end{equation}
Both terms in the sum and integral above are positive, and thus the
integrals of the sums of the individual terms are also bounded.

We can use the inequality 
\begin{equation*}
  \abs{R}^\alpha + \abs{S}^\alpha \ge C_\alpha
  \abs{\left(\abs{R}-\abs{S}\right)}^\alpha=C_\alpha\abs{R-S}^\alpha 
\end{equation*}
for some constant $C_\alpha$ depending on $\alpha$, to get the bound
\begin{equation}
  \label{eq:highbnd1}
  \int_0^T \Dx\sum_j \chi_j \tc_j \abs{R_j-S_j}^{2+\alpha}\,dt \le
  C_{\alpha,T,a,b}. 
\end{equation}
Since $C_1<c(u)<C_2$ (cf.~\eqref{c_est}) we have proved the following lemma.
\begin{lemma}
  \label{lem:1alphabnd}
  Let $\alpha\in [0,1)$, and assume that $R_0$ and $S_0$ are
  nonpositive, and in $L^1(\R)\cap L^2(\R)$. Then 
  we have the estimate 
  \begin{equation}
    \label{eq:dualphaplus}
    \int_0^T \Dx \sum_{j=j_a}^{j_b} c'(u_j^+) \abs{\Dp \uj{-}}^{2+\alpha}\,dt \le
    C_{\alpha,T,a,b},
  \end{equation}
  where $j_a\Dx\in[a-1,a)$ and $j_b\Dx \in(b,b+1]$.
\end{lemma}

\medskip
{\bf Acknowledgements.} The authors gratefully acknowledge the  
hospitality of the Mittag-Leffler Institute, Sweden, creating a great  
working environment for research, during the Fall of 2005.


\end{document}